\documentclass[centertags,leqno,11pt]{article}

\usepackage{amsmath}
\usepackage{amssymb}
\usepackage{latexsym}
\usepackage{amsxtra}
\usepackage{amscd}
\usepackage{theorem}
\usepackage[all]{xy}
\input xy
\xyoption{all}
\entrymodifiers={+!!<0pt,\fontdimen22\textfont2>}

\usepackage{graphics}
\usepackage{epic}

\usepackage{a4wide, mathrsfs, amsfonts, latexsym, dsfont}

\theoremstyle{change}

{\theorembodyfont{\slshape}

\newtheorem{thm}{Theorem.}[section]
\newtheorem{cor}[thm]{Corollary.}
\newtheorem{lem}[thm]{Lemma.}
\newtheorem{prop}[thm]{Proposition.}

}

{\theorembodyfont{\rmfamily}

\newtheorem{rem}[thm]{Remark.}

}

\renewcommand{\em}{\sl}

\newcommand{\proof}{\noindent {\bf Proof:\ }}
\newcommand{\Endproof}{\hspace*{\fill} $\Box$ \vspace{1ex} \noindent }

\makeatletter
\renewcommand{\subsection}{\@startsection{subsection}{2}%
{\z@}{-3.25ex plus -1ex minus-.2ex}{-1em}{\bf}} \makeatother

\newcommand{\ZZ}{\mathbb{Z}}

\newcommand{\FF}{\mathbb{F}}
\renewcommand{\AA}{\mathbb{A}}

\newcommand{\W}{\mathcal{W}}

\newcommand{\GL}{{\rm GL}}

\newcommand{\Aut}{{\rm Aut}}

\numberwithin{equation}{section}
\numberwithin{thm}{section}
\theoremstyle{plain}

\title{Rigid $G_2$-Representations and motives of Type $G_2$}

\author{Michael Dettweiler\footnote{The first author gratefully acknowledges 
financial support from the DFG-Heisenberg Grant DE-1442}, Johannes Schmidt\footnote{The second author gratefully acknowledges 
financial support from the DFG-Forschergruppe 1920 "Symmetrie, Geometrie und Arithmetik", Heidelberg--Darmstadt}}

\begin{document}

\maketitle

\section{Introduction}

\subsection*{Background.}
Pure motives in the sense of Grothendieck form a Tannaka category if the standard conjectures on algebraic cycles hold 
(cf.~\cite{Andre04}).
The category of motives for absolute Hodge-cycles \cite{DelMil82} and  the category of motives for motivated cycles \cite{Andre96}
each form unconditionally a Tannaka category.  By the Tannakian formalism, 
a pure motive $M$ admits then a motivic Galois group $G_M,$ which measures how the tensor 
constructions of $M$ decompose (assuming the standard conjectures, when working with 
Grothendieck motives).
In \cite{Serre}~8.8, Serre posed the question
\begin{quote}
 Exist-t-il un motif $M$ tel que $G_{M}$ soit un groupe simple de type exceptionnel $G_2$ (ou $E_8$)?                              
\end{quote}
The $G_2$-case has been positively answered in \cite{Dettweiler1} as follows (cf.\ \cite{Yun14} for the $E_8$-case):
All rigid local systems on $\mathbb{P}_{\bar{\mathbb{Q}}}^1\setminus \{ 0,1,\infty \}$ of rank $7$ with monodromy group 
$G_2$ where classified up to isomorphisms (see \cite{Dettweiler1} Thm.\ 1.3.1).
There are exactly $5$ series of isomorphism classes of such rigid local systems, distinguished by the 
 Jordan decomposition of the  unipotent part of the local monodromy at $\infty$ (the series being parametrized by the eigenvalues appearing there). 

 The motivic description of rigid local systems 
 in \cite{Katz} Chapter 8 describes these local systems as geometric generic fibre of a certain pure part of a direct image to $\mathbb{A}_R^1 \setminus \{0,1\}$ in $\ell$-adic cohomology, where $R/\mathbb{Z}$ is generically generated by the eigenvalues appearing in the geometric monodromy (see \cite{Dettweiler1} Thm.\ 2.4.1).
Using this motivic description of a member $\mathcal{V}$ of the isomorphism class whose monodromy at $\infty$ is given by the unipotent Jordan block of length $7$, a family of rational motives for motivated cycles $\{M_s\}_{s\in \mathbb{A}^1\setminus \{ 0,1 \}}$ over $\mathbb{A}^1\setminus \{ 0,1 \}$ was constructed s.t.\ the $\ell$-adic realization of $M_s$ is canonically isomorphic to the stalk $\bar{s}^*\mathcal{V}(3)$ of the Tate twist over $s$.
The relation between the motivic Galois groups of such a family of motives and the action on the Betti realizations of the fundamental group of the base space is described in \cite{Andre96} Thm.\ 5.2.
Using quite general arguments around this description allowed Dettweiler and Reiter to conclude that the motivic Galois Group of $M_s$ is contained in the monodromy Group $G_2$ of $\mathcal{V}(3)$ (and agrees with this monodromy group for $s$ outside a meagre subset of the base space).
\\
Further, the Galois-action on the $\ell$-adic realization of $M_s$ was studied by Dettweiler and Katz:
Using the indecomposable action of the geometric monodromy at $\infty$, the irreducibility of the Galois-action could be shown.
This in turn was the key to an explicit description of infinitely many $s\in \mathbb{Q}\setminus \{0,1\}$ s.t.\ the motivic Galois group of $M_s$  agrees with $G_2$.

In contrast, the other four series of isomorphism classes of rigid local systems with monodromy $G_2$ lack an indecomposable monodromy action at $\infty$.
Further, their motivic descriptions lead only to families of motives for motivated cycles $\{M_s\}_s$ over the number fields generated by the eigenvalues of the respective monodromies, which are all bigger than $\mathbb{Q}$.

\subsection*{Aim and results of the paper.}
In this paper we want to study the case of a rigid local system $\mathcal{V}$ of type $G_2$ whose monodromy at $\infty$ is given by
\begin{equation*}
 {\rm diag}(\zeta_7,\zeta_7^2,\zeta_7^3,\zeta_7^4,\zeta_7^5,\zeta_7^6,1), \end{equation*}
  cf.~\cite{Dettweiler1}
 Thm.\ 1.3.2.
It follows that the corresponding family of motives $\{ M_s \}_s$ is  defined over $\mathbb{Q}(\zeta_7)$ and has  coefficients 
also in $\mathbb{Q}(\zeta_7)$.

Apart from minor differences, the arguments that $G_2$ is an upper bound for the motivic Galois group of $M_s$ apply in our case, as well.
However, the discussion leading to the irreducibility of the Galois action on the $\ell$-adic realization fails completely due to the decomposability of the geometric monodromy action at $\infty$.
Our way around this is to show that the mod-$\ell$ reduction $\bar{\mathcal{V}}_{\mathbb{Q}(\zeta_7)}$ of the extension $\mathcal{V}_{\mathbb{Q}(\zeta_7)}$ to the base field $\mathbb{Q}(\zeta_7)$ (given by the motivic description of $\mathcal{V}$) has irreducible stalks.
To see this, we need the mod-$\ell$ reduction of $\mathcal{V}$ to descend even to a locally constant sheaf on $\mathbb{A}_{\mathbb{Q}}^1 \setminus \{ 0,1\}$ with monodromy group $G_2(\mathbb{F}_\ell)$.
The reason for this is that we want the Frobenius lift of a place over which the corresponding Galois representation is unramified to act transitive on the set of eigenspaces for eigenvalues $\neq 1$ of the geometric monodromy at $\infty$.
Using the rigidity method of inverse Galois theory (e.g.\ see \cite{Matzat} Chapter I), a unique such descent (or extension in the language of Galois theory) $\bar{\mathcal{V}}_{\mathbb{Q}}$ with monodromy group $G_2(\mathbb{F}_\ell)$ is provided in Prop.\ \ref{prop: descent to Q}. 
Similar to \cite{Dettweiler1} Prop.\ A.4, a careful study of the restrictions to inertia subgroups of the $G_{\mathbb{Q}}$-action on the stalks over the infinitely many $s \in \mathbb{Q} \setminus \{ 0,1 \}$ (cf.~Rem.~\ref{rem: infinitely many s}) satisfying
\begin{center}
 (\dag)~ \parbox[t]{12cm}{there are prime numbers $p, q \neq 2, 7$ with $p \equiv 3 ~{\rm or}~ 5 ~{\rm mod}~ 7$ such that $\nu_p(s) < 0 < \nu_q(s - 1)$ and $7 \nmid \nu_p(s)$}
\end{center}
then leads us to a Hilbert-irreducibility result for this action (for the exact statement, see Thm.\ \ref{thm: hilbert irreducibility} below):

\phantom{leere Zeile}

\noindent {\bf Theorem~I.} {\it
For the infinitely many rational numbers $s \neq 0,1$ satisfying {\rm(\dag)} and a prime $\ell \gg 0$,
the image of the Galois-action on the stalk $\bar{s}^*\bar{\mathcal{V}}_{\mathbb{Q}}$ is again $G_2(\mathbb{F}_\ell)$. 
}

\phantom{leere Zeile}

Now $G_2(\mathbb{F}_\ell)$ is a simple group acting irreducibly on the stalks of $\bar{\mathcal{V}}_{\mathbb{Q}}$, so the Galois action on corresponding stalks of $\bar{\mathcal{V}}_{\mathbb{Q}(\zeta_7)}$ is still irreducible, too.
This in turn implies the desired irreducibility of the Galois action on the $\ell$-adic realization of $M_s$.
From this we conclude that the image of this Galois-action up to a twist of a character $\varepsilon$ of order at most $2$ is in fact $G_2$.
Let $k / \mathbb{Q}(\zeta_7)$ be the extension of degree at most $2$ trivializing $\varepsilon$, i.e., the image of the Galois-action on the $\ell$-adic realization of $M_s$ over $k$ is $G_2$ (for the exact statement, see Thm.\ \ref{thm: tate realization} below):

\phantom{leere Zeile}

\noindent {\bf Theorem~II.} {\it
For the infinitely many rational numbers $s \neq 0,1$ satisfying {\rm(\dag)} and a prime $\ell \gg 0$, the absolute Galois group $G_k$ acts as a Zariski dense subgroup of $G_2(\bar{\mathbb{Q}}_{\ell})$ on the $\ell$-adic realization $H_\ell^\bullet(M_s)$ of the motive for motivated cycles $M_s$ over $k$.
}

\phantom{leere Zeile}

The Tate conjecture implies that this image is the motivic Galois group itself, so a relatively self contained r\'{e}sum\'{e} is (see Cor.\ \ref{cor: numerical motives under tate conjecture} below):

\phantom{leere Zeile}

\noindent {\bf Corollary.} {\it
Suppose the Tate conjecture holds.
For the infinitely many rational numbers $s \neq 0,1$ satisfying {\rm(\dag)}, the motivic Galois group of the numerical motive $M_s/k$ is of type $G_2$. 
}

\phantom{leere Zeile}

A similar result holds unconditionally for the category of motives for absolute Hodge cycles, resp. for the category 
of motives for motivated cycles (Thm.~ \ref{thm: motives for absolute hodge cycles}).

\subsection*{Acknowledgements.}
We would like to thank Armin Holschbach, Stefan Reiter and Jakob Stix for helpful discussions or comments on earlier versions of this paper.

\section{Construction of sheaves with monodromy $G_2$}\label{con}

\subsection{Introduction of the basic objects.}\label{basicob}
Fix a prime number $\ell$.
For $N$ a natural number denote by $R_N$ the ring $\mathbb{Z}[\zeta_7, \frac{1}{14 \cdot \ell \cdot N}]$.
We follow \cite{Dettweiler1} to construct in a geometric (motivic) way a lisse sheaf on $\mathbb{A}_{R_N}^1 \setminus \{ 0,1 \}$ of rank $7$ whose monodromy representation has a Zariski dense image inside $G_2(\bar{\mathbb{Q}}_{\ell})$ embedded into ${\rm GL}_7(\bar{\mathbb{Q}}_\ell)$ as in \cite{Aschbacher}.

Set 
\begin{equation*}
\Delta := \prod_{i,j}(X_i-j) \prod_k(X_{k+1}-X_k),
\end{equation*}
where $i$ runs through $\lbrace 1,...,7 \rbrace$, $j$ through $\lbrace 0,1 \rbrace$ and $k$ runs through $\lbrace 1,...,6 \rbrace$.
Denote by ${\rm Hyp}_N$ the hypersurface in $\mathbb{G}_{{\rm m},R_N} \otimes_{R_N} (\mathbb{A}_{R_N}^7 \setminus v(\Delta))$ defined by
\begin{equation*}
Y^{14} = \prod_{i, j}(X_i - j)^{e(i,j)} \prod_{k}(X_{k+1} -X_k)^{f(k)}
\end{equation*}
($i,j,k$ as above), where $Y$ denotes the coordinate function of $\mathbb{G}_{{\rm m},R_N}$, the $X_i$'s denote the coordinate functions of $\mathbb{A}_{R_N}^7$ and the exponents are as follows:

\begin{displaymath}
\begin{array}{|l|l|l|l|l|l|l|}
\hline
e(1,0) & e(2,0) & e(3,0) & e(4,0) & e(5,0) & e(6,0) & e(7,0)
\\
\hline
7 & 0 & 7 & 0 & 7 & 0 & 7
\\
\hline \hline
e(1,1) & e(2,1) & e(3,1) & e(4,1) & e(5,1) & e(6,1) & e(7,1)
\\
\hline
13 & 5 & 0 & 5 & 0 & 5 & 0
\\
\hline
\end{array}
\end{displaymath}
and
\begin{displaymath}
\begin{array}{|l|l|l|l|l|l|}
\hline
f(1) & f(2) & f(3) & f(4) & f(5) & f(6)
\\
\hline
11 & 3 & 1 & 13 & 5 & 9
\\
\hline
\end{array}.
\end{displaymath}
Let
\begin{equation*}
\xymatrix{
\pi_N:{\rm Hyp}_N \ar[r] & \mathbb{A}_{R_N}^1 \setminus \lbrace 0,1 \rbrace
}
\end{equation*}
be the morphism induced by the inclusion
\begin{equation*}
\xymatrix{
R_N[X_7,\frac{1}{X_7(X_7 - 1)}] \ar@{^{(}->}[r] & R_N[X_1,...,X_7,Y,\frac{1}{Y \Delta}].
}
\end{equation*}
Denote by $\boldsymbol{\mu}_{14}$ the group of $14^{th}$ roots of unity in $R_N$.
It acts on ${\rm Hyp}_N$ via $Y \mapsto \omega Y$ (for $\omega \in \boldsymbol{\mu}_{14}$), inducing an $\boldsymbol{\mu}_{14}$-action on $\mathds{R}^q\pi_{N,!}\bar{\mathbb{Q}}_{\ell}$.
Fix an embedding of $R_N$ into $\bar{\mathbb{Q}}_{\ell}$.
Let $\chi$ be the $\bar{\mathbb{Q}}_{\ell}$-valued character of $\boldsymbol{\mu}_{14}$ given by this embedding.
By $(\mathds{R}^q\pi_{N,!}\bar{\mathbb{Q}}_{\ell})^{\chi}$ we mean the $\chi$-eigenspace of the $\boldsymbol{\mu}_{14}$-action.
Then by \cite{Dettweiler1} Thm.\ 2.3.1 and Thm.\ 2.4.1, $(\mathds{R}^q\pi_{N,!}\bar{\mathbb{Q}}_{\ell})^{\chi}$ is a lisse $\bar{\mathbb{Q}}_{\ell}$-sheaf on $\mathbb{A}_{R_N}^1 \setminus \lbrace 0,1 \rbrace$ of mixed weights $\leq q$ and
\begin{equation*}
\mathcal{V} := W_6[(\mathds{R}^6\pi_{N,!}\bar{\mathbb{Q}}_{\ell})^{\chi}]_{\mid \mathbb{A}_{\bar{\mathbb{Q}}}^1 \setminus \{ 0,1 \}}
\end{equation*}
(here $W_\bullet$ denotes the weight filtration) is cohomologically rigid (see \cite{Katz} (5.0.1)) of rank $7$ with induced monodromy representation $\rho_{\mathcal{V}}$ Zariski dense in $G_2(\bar{\mathbb{Q}}_{\ell})$ embedded into ${\rm GL}_7(\bar{\mathbb{Q}}_\ell)$ as in \cite{Aschbacher}.
Further, its monodromy tuple $(g_0,g_1,g_{\infty})$ has the following Jordan canonical forms:
\begin{eqnarray}
{\rm JCF}(g_0) &=& {\rm diag}(-1,-1,-1,-1,1,1,1), \nonumber
\\
{\rm JCF}(g_1) &=& 
{\bf U}(2)\oplus{\bf U}(2)\oplus{\bf U}(3) \nonumber
\\
{\rm JCF}(g_{\infty}) &=& {\rm diag}(\zeta_7,\zeta_7^2,\zeta_7^3,\zeta_7^4,\zeta_7^5,\zeta_7^6,1), \nonumber
\end{eqnarray}
where ${\bf U}(k)$ denotes the unipotent Jordan block of length $k$.

Throughout the paper we use the following notation: If $\W$ is a lisse or \'{e}tale locally constant sheaf on 
a connected scheme $X,$ then $\rho_\W:\pi_1(X, \bar{s})\to \Aut(\W_{\bar{s}})$ denotes the monodromy representation 
of $\W$ with base point $\bar{s}$.  
Set
\begin{equation*}
 \mathcal{V}_k := W_6[(\mathds{R}^6\pi_{N,!}\bar{\mathbb{Q}}_{\ell})^{\chi}]_{\mid \mathbb{A}_{k}^1 \setminus \{ 0,1 \}}
\end{equation*}
for $k$ a number field containing $\zeta_7$ (similar define $\mathcal{V}_{R_N}, \dots$).
By definition, this is a descent of $\mathcal{V}$ to $\mathbb{A}_{k}^1 \setminus \{ 0,1 \}$.
Since the geometric fundamental group is a normal subgroup of $\pi_1^{{\rm \acute{e}t}}(\mathbb{A}_{k}^1 \setminus \{ 0,1 \})$, the image of the induced representation $\rho_{\mathcal{V}_k}$ normalizes the image of $\rho_{\mathcal{V}}$.

By construction, $\rho_{\mathcal{V}_{R_N}}$ is  ${\rm GL}_7(\bar{\mathbb{Q}}_{\ell})$-conjugate to a representation with values in ${\rm GL}_7(\mathbb{Z}_{\ell}(\zeta_7))$.
Thus we may assume that $\rho_{\mathcal{V}}, \rho_{\mathcal{V}_{R_N}}, \rho_{\mathcal{V}_{k}}, \dots$ have values in ${\rm GL}_7(\mathbb{Z}_{\ell}(\zeta_7))$.
Define the mod $\ell$ reductions $\bar{\mathcal{V}}, \bar{\mathcal{V}}_{R_N}, \bar{\mathcal{V}}_{k}, \dots$ as the \'{e}tale locally constant sheaves given by the naive mod $\ell$ reductions of the above representations.
By construction, the monodromy representation of $\bar{\mathcal{V}}$ has image contained in $G_2(\bar{\mathbb{F}}_\ell)$ embedded into ${\rm GL}_7(\bar{\mathbb{F}}_\ell)$ as in \cite{Aschbacher}.
In the following two subsections we will see that the monodromy group is $G_2(\mathbb{F}_\ell)$.

\subsection{Computation of the monodromy of $\bar{\mathcal{V}}$.}\label{sect: computation of monodromy}
If $G$ is an algebraic group defined over $\ZZ$ (e.g. $G=\GL_n$ or $G_2\leq \GL_7$)
and if $\ell^n$ is a power of a prime number $\ell,$ we write $G(\ell^n)$ for  the group
$G(\FF_{\ell^n}).$
From now on let $\ell \neq 2,3,7$.
Let $\boldsymbol{\sigma} := (\bar{g}_0, \bar{g}_1, \bar{g}_{\infty})$ be the monodromy tuple of $\bar{\mathcal{V}}$.
Then the Jordan canonical forms of the tuple $\boldsymbol{\sigma}$ remain untouched:

\begin{lem}\label{lem: jordan forms remain untouched}
Let ${\ell} \neq 2,3,7$ and $i=0,1,\infty$.
Then ${\rm JCF}(\bar{g}_i)$ is the mod-$\ell$-reduction of ${\rm JCF}(g_i)$.
\end{lem}

\proof
We start with $\bar{g}_0$.
It is an involution so its Jordan canonical form is determined by its trace, which is the mod $\ell$ reduction of ${\rm tr}(g_0) = -1$.
The other possible traces of an involution of rank $7$ are $1,\pm2,\pm3,\pm4,\pm5$ and $-7$.
Since $\ell \neq 2,3$, the traces of the mod-$\ell$-reductions of the other types of involutions  are all different from $-1$, i.e. ${\rm JCF}(\bar{g}_0) = \overline{{\rm JCF}(g_0)}$ follows.\\
By looking at the minimal polynomial it is clear that $\bar{g}_1$ is unipotent with Jordan blocks of length less than three.
If it has more than three blocks, the dimension of the eigenspace ${\rm Eig}(\bar{g}_1,1)$ would be at least $4$.
Now ${\rm Eig}(\bar{g}_0,-1)$ has dimension $4$.
In particular the intersection of both eigenspaces is non trivial and thus the product relation $\bar{g}_0 \bar{g}_1 \bar{g}_{\infty} = \textbf{1}$ implies that $\bar{g}_{\infty}$ has $-1$ as eigenvalue which clearly is not the case.
The only remaining possibility except ${\rm JCF}(\bar{g}_1) = \overline{{\rm JCF}(g_1)}$ is two blocks of length three and one of length one.
But then ${\rm rk}((\bar{g}_1 - \textbf{1})^2)$ would be greater then ${\rm rk}((g_1 - \textbf{1})^2)$ which is impossible.\\
Finally, $\ell \neq 7$, so the characteristic polynomial of $\bar{g}_{\infty}$ is still separable and ${\rm JCF}(\bar{g}_{\infty})$ is the reduction of $\overline{{\rm JCF}(g_{\infty})}$, as well.
\Endproof

Recall that $\boldsymbol{\sigma}$ is linearly rigid if and only if each triple of elements in the ${\rm GL}_7(\bar{\mathbb{F}}_{\ell})$-conjugacy classes of the $\bar{g}_i$'s that satisfies the above product relation is even simultaneously conjugate to $\boldsymbol{\sigma}$.
Using Lem.\ \ref{lem: jordan forms remain untouched} we get:

\begin{lem}\label{lem: linear rigidity}
The tuple $\boldsymbol{\sigma}$ is absolutely irreducible and linearly rigid.
\end{lem}

\proof
Recall that the Katz algorithm for rigids \cite{Katz} Chapter 6 works over arbitrary algebraically closed fields, too (see \cite{Dettweiler3} Sect.\ 4).
By Lem.\ \ref{lem: jordan forms remain untouched} we can apply it simultaneously to $(g_0,g_1,g_\infty)$ and $\boldsymbol{\sigma}$.
Since $\mathcal{V}$ was constructed via this algorithm ``applied backwards'', it stops with an absolutely irreducible and cohomologically rigid tuple of rank $1$ in the former case.
Thus, the same is true in the latter case, which completes the proof.
\Endproof

\begin{rem}\label{rem: maximal subgroups}
For $n \gg 0$ the image of $\rho_{\bar{\mathcal{V}}}$ lies in $G_2(\ell^n)$ embedded into ${\rm GL}_7(\ell^n)$ as in \cite{Aschbacher}.
We write $V(\ell^n)$ for the corresponding $\mathbb{F}_{\ell^n}$-structure of a stalk $\bar{s}^*\bar{\mathcal{V}}$, i.e., $\rho_{\bar{\mathcal{V}}}$ factors over an image of $G_2(\ell^n)$ inside ${\rm GL}(V(\ell^n))$.
\\
In \cite{Kleidman}, Kleidman has shown that every maximal subgroup $M$ of $G_2({\ell}^n)$ for ${\ell} \neq 2,3$ is $G_2({\ell}^n)$-conjugate to one of the groups contained in the following list (in Kleidman's notation, our vector space $V(\ell^n)$ would be the subspace $W$ of the octonions $V$ - note that this is exactly the space where the Dickson form lives on used by Aschbacher):
\begin{center}
\begin{tabular}{|l|l|}
\hline
\textbf{Group:} & \textbf{Remarks:}
\\
\hline \hline
$P_{\alpha}, P_{\beta}$ & maximal parabolic subgroups
\\
\hline
$C_{G_2({\ell}^n)}(\iota)$ & centralizer of the involution $\iota$
\\
\hline
$K_{\varepsilon} \cong L_{\varepsilon} \rtimes \mathbb{Z} / 2 $ & $\varepsilon = \pm, ~ L_+ \cong {\rm SL}_3({\ell}^n), ~ L_- \cong {\rm SU}_3({\ell}^n)$
\\
\hline
${\rm PGL}_2({\ell}^n)$ & for ${\ell} \geq 7, {\ell}^n \geq 11$
\\
\hline
$ (\mathbb{Z} / 2)^3 \dot{~} L_3(2)$ & for $n = 1$
\\
\hline 
$L_2(8)$ & for ${\ell} \geq 5, \mathbb{F}_{{\ell}^n} = \mathbb{F}_{{\ell}}(\omega)$ with $\omega^3 - 3\omega + 1 = 0$
\\
\hline 
$L_2(13)$ & for ${\ell} \neq 13, \mathbb{F}_{{\ell}^n} = \mathbb{F}_{{\ell}}(\sqrt{13})$
\\
\hline 
$G_2(2)$ & for ${{\ell}^n}={\ell} \geq 5$
\\
\hline 
$J_1$ & for ${{\ell}^n}=11$
\\
\hline 
$C_{G_2({\ell}^n)}(\phi_{\ell^m}) = G_2({\ell}^m)$ & for $\alpha = \frac{n}{m}$ prime
\\
\hline 
\end{tabular}
\end{center}
Here we denote by $ (\mathbb{Z} / 2)^3 \dot{~} L_3(2)$ a suitable non-split group extension of $(\mathbb{Z} / 2)^3$ with $L_3(2)$.
We do not need an exact definition of $P_{\alpha}$ and $P_{\beta}$ at this point, but it can be found later in section \ref{image}.
\end{rem}

We will use this list to prove the following proposition:

\begin{prop}\label{prop: image G2}
Let ${\ell} \neq 2,3,7,11$ and $13$ and suppose ${\rm im}({\rho}_{\bar{\mathcal{V}}}) \leq G_2({\ell}^n)$.
Then ${\rho}_{\bar{\mathcal{V}}}$ has image $G_2({\ell}^n)$ or ${\rm im}({\rho}_{\bar{\mathcal{V}}})$ is contained in a subgroup $G_2({\ell}^n)$-conjugate to $C_{G_2({\ell}^n)}(\phi_{\ell^m}) = G_2({\ell}^m)$ for $\alpha = \frac{n}{m}$ a prime number.
\end{prop}

\proof
We simply have to check for every type of group in the list, whether the group contains ${\rm im}({\rho}_{\bar{\mathcal{V}}})$.
\\
\textbf{(i)}
By Lem.\ \ref{lem: linear rigidity}, ${\rm im}({\rho}_{\bar{\mathcal{V}}})$ acts irreducibly on $V({\ell}^n)$.
The first three groups in the table act reducibly, so ${\rm im}({\rho}_{\bar{\mathcal{V}}})$ is not contained in any subgroup $G_2({\ell}^n)$-conjugate to one of these subgroups.
\\
\textbf{(ii)}
As in the proof of \cite{Kleidman} Lem.\ 4.2.1 $V({\ell}^n)$ as ${\rm PGL}_2({\ell}^n)$-module is isomorphic to the $\mathbb{F}_{{\ell}^n}$-vector space of homogeneous polynomials in $X$ and $Y$ of degree 6 as representation of
\begin{equation*}
\langle \rho({\rm SL}_2({\ell}^n)), \mu^{-3} \rho({\rm diag}(\mu,1)) \rangle, 
\end{equation*}
where $\mu$ is a non square in $\mathbb{F}_{{\ell}^n}$ and $\rho$ is the ${\rm GL}_2({\ell}^n)$ representation given by
\begin{equation}\label{eq: PSL-action}
\begin{pmatrix}
a & b \\ c & d
\end{pmatrix}
: X^i Y^j \mapsto (aX + bY)^i (cX + dY)^j.
\end{equation}
The non trivial unipotent elements in ${\rm GL}_2(\ell^n)$ form a single conjugacy class.
Hence the Jordan canonical form of each unipotent element in our ${\rm PGL}_2({\ell}^n)$-representation consists of a single Jordan block of length $7$.
Since this is not the case for $\bar{g}_1$, no maximal subgroup $G_2({\ell}^n)$-conjugate to ${\rm PGL}_2({\ell}^n)$ contains ${\rm im}({\rho}_{\bar{\mathcal{V}}})$.
\\
\textbf{(iii)}
Finally, we can exclude all other maximal subgroups but $C_{G_2({\ell}^n)}(\phi_{\alpha})$ in the above table by comparing their orders with our assumptions on $\ell$:
From \cite{Atlas} we get the orders of the groups $ (\mathbb{Z} / 2)^3 ~ \! ^{\cdot}  L_3(2)$, $L_2(8)$, $L_2(13)$ and $J_1$.
The order of $G_2(2)$ can be computed using \cite{Carter} Thm.\ 9.4.10, Cor.\ 10.2.4 and Prop.\ 10.2.5.
These information together give us the following list.
\begin{center}
\begin{tabular}{|l|l|}
\hline
\textbf{Group:} & \textbf{Order:}
\\
\hline \hline
$ (\mathbb{Z} / 2)^3 ~ \! ^{\cdot}  L_3(2)$ & $2^6 \cdot 3 \cdot 7$
\\
\hline
$L_2(8)$ & $2^3 \cdot 3^2 \cdot 7$
\\
\hline
$L_2(13)$ & $2^2 \cdot 3 \cdot 7 \cdot 13$
\\
\hline
$J_1$ & $2^3 \cdot 3 \cdot 5 \cdot 7 \cdot 11 \cdot 19$
\\
\hline
$G_2(2)$ & $2^6 \cdot 3^3 \cdot 7$
\\
\hline
\end{tabular}
\end{center}
Since $J_1$ appears only on the upper list of maximal subgroups for ${\ell} = 11$, non of these groups contain an element of order $\ell$ for ${\ell} \neq 2,3,7,11 ~ {\rm and} ~ 13$.
But $\bar{g}_1$ has order ${\ell}$ by Lem.\ \ref{lem: jordan forms remain untouched}, which completes the proof.
\Endproof

Iterative application of this gives us an $n$ s.t.\ ${\rm im}(\rho_{\bar{\mathcal{V}}})$ is $G_2(\ell^n)$.
In particular, our monodromy tuple $\boldsymbol{\sigma}$ generates $G_2(\ell^n)$.
We will see in Cor.\ \ref{cor: n=1} below that $n$ is actually just $1$.

\subsection{Descent of $\bar{\mathcal{V}}$ to $\mathbb{A}_{\mathbb{Q}}^1 \setminus \{ 0,1 \}$.}
We turn our attention again to the locally constant sheaf $\bar{\mathcal{V}}$ on $\mathbb{A}_{\bar{\mathbb{Q}}}^1 \setminus \{ 0,1 \}$.
In this subsection we will prove:

\begin{prop}\label{prop: descent to Q}
Let ${\ell} \neq 2,3,7,11$ and $13$.
Then $\bar{\mathcal{V}}$ descends to a local system $\bar{\mathcal{V}}_{\mathbb{Q}}$ on $\mathbb{A}_{\mathbb{Q}}^1 \setminus \{ 0,1 \}$, satisfying
\begin{equation*}
{\rm im}(\rho_{\bar{\mathcal{V}}_{\mathbb{Q}}}) = {\rm im}(\rho_{\bar{\mathcal{V}}}) = G_2(\ell) .
\end{equation*}
\end{prop}

\begin{rem}\label{rem: local system with coefficients in the prime field}
Until now, we treated $\bar{\mathcal{V}}$ as a local system of $\bar{\mathbb{F}}_\ell$-vector spaces that is already defined over $\mathbb{F}_{\ell^n}$ (with $n = 1$: see Cor.\ \ref{cor: n=1}, below). 
Yet, we will treat $\bar{\mathcal{V}}_{\mathbb{Q}}$ as a local system of $\mathbb{F}_\ell$-vector spaces.
So strictly speaking, $\bar{\mathcal{V}}_{\mathbb{Q}} \otimes \bar{\mathbb{F}}_\ell$ is the descent of $\bar{\mathcal{V}}$ to $\mathbb{A}_{\mathbb{Q}}^1 \setminus \{ 0,1 \}$.
\end{rem}

\noindent \textbf{Proof of Prop.\ \ref{prop: descent to Q}:}
Let us first recall some facts from Galois theory:
Let $\mathbb{S}$ be the set of places of $\bar{\mathbb{Q}}(t)$ corresponding to the $t$- and $(t-1)$-adic valuation and the valuation at infinity.
As in \cite{Matzat} I Thm.\ 2.2 denote by $\bar{M}_{\mathbb{S}}$ the maximal algebraic extension field of $\bar{\mathbb{Q}}(t)$ inside $\overline{\mathbb{Q}(t)}$ unramified outside $\mathbb{S}$.
Then ${\rm Gal}(\bar{M}_{\mathbb{S}} / k(t))$ is canonically isomorphic to $\pi_1^{{\rm \acute{e}t}}(\mathbb{A}_{k}^1 \setminus \{ 0,1 \}, \bar{0})$ for $k$ an algebraic extension field of $\mathbb{Q}$, where $\bar{0}$ denotes the geometric point corresponding to the choice of an embedding of $k(t)$ to $\overline{\mathbb{Q}(t)}$.
In particular $\rho_{\bar{\mathcal{V}}}$ gives rise to an intermediate field of $\bar{M}_{\mathbb{S}} / \bar{\mathbb{Q}}(t)$, Galois over $\bar{\mathbb{Q}}(t)$ with group isomorphic to $G_2(\ell^n)$ ($n$ as in the last subsection).
As in the Hurwitz classification of such intermediate fields in \cite{Matzat} I.4.1 denote this field by $\bar{N}_{\boldsymbol{\sigma}}$ (i.e., $\rho_{\bar{\mathcal{V}}}$ corresponds to $\psi_{\boldsymbol{\sigma}}$ in loc.\ cit.).
Now, if $\mathbb{Q}(t)$ is a field of definition of $\bar{N}_{\boldsymbol{\sigma}} /_{G_2(\ell^n)} \bar{\mathbb{Q}}(t)$ (see \cite{Matzat} I.3.1), we get an epimorphism
\begin{equation*}
 \rho_{\bar{\mathcal{V}}_{\mathbb{Q}}}:
 \xymatrix{
  {\rm Gal}(\bar{M}_{\mathbb{S}} / \mathbb{Q}(t)) \ar@{->>}[r] &
  G_2(\ell^n)
 }
\end{equation*}
which, after a possible modification by an automorphism of $G_2(\ell)$, restricts to $\rho_{\bar{\mathcal{V}}}$.
As a result the locally constant sheaf $\bar{\mathcal{V}}_{\mathbb{Q}}$ on $\mathbb{A}_{\mathbb{Q}}^1 \setminus \{ 0,1 \}$ defined by $\rho_{\bar{\mathcal{V}}_{\mathbb{Q}}}$ serves as the desired descent of $\bar{\mathcal{V}}$.
\\
Thus, it remains to prove that $\mathbb{Q}(t)$ is in fact a field of definition of $\bar{N}_{\boldsymbol{\sigma}} /_{G_2(\ell^n)} \bar{\mathbb{Q}}(t)$:
According to Lem.\ \ref{lem: matzat rigid} below, $\boldsymbol{\sigma}$ is rigid in the sense of \cite{Matzat}, i.e., for each triple $\{h_i\}$ in $G_2(\ell^n)^3$ satisfying the product relation $\bar{g}_0^{h_0} \bar{g}_1^{h_1} \bar{g}_\infty^{h_\infty} = \boldsymbol{1}$, there is even an element $h$ of $G_2(\ell^n)$ with $\bar{g}_i^h = \bar{g}_i^{h_i}$ simultaneously for all $i$.
Further, according to Lem.\ \ref{lem: rational} below, $\boldsymbol{\sigma}$ is also rational in the sense of \cite{Matzat}, i.e., for all primes $p$ not dividing the order of $G_2(\ell^n)$ the $p$ power map preserves the $G_2(\ell^n)$-conjugacy class of each $\bar{g}_i$.
As a corollary to the arguments leading to the rigidity of $\boldsymbol{\sigma}$ (see Rem.\ \ref{rem: normalizer2}) we will also see that actually $n$ is just $1$ (see Cor.\ \ref{cor: n=1}, below). 
Then $\mathbb{Q}(t)$ is in fact a field of definition of $\bar{N}_{\boldsymbol{\sigma}} /_{G_2(\ell)} \bar{\mathbb{Q}}(t)$ by the Rigidity Theorem (see \cite{Matzat} I Thm.\ 4.8), which completes the proof.
\Endproof

We start with checking the rationality of $\boldsymbol{\sigma}$:

\begin{lem}\label{lem: rational}
The tuple $\boldsymbol{\sigma}$ is rational with respect to $G_2(\ell^n)$ in the sense of \cite{Matzat}.
\end{lem}

\proof
By \cite{Matzat} I Prop.\ 4.4 this is equivalent to the statement that every complex irreducible character of $G_2(\ell^n)$ has rational values for $\bar{g}_0, \bar{g}_1$ and $\bar{g}_{\infty}$.
But this is straightforward, using the character table in \cite{changree}.
Thus $\boldsymbol{\sigma}$ is indeed rational.
\Endproof

For the rigidity, we need a little preparation:

\begin{rem}\label{rem: h infty p}
Let $p \equiv 3 ~{\rm or}~ 5 ~{\rm mod}~7$ be a prime number not dividing the order of $G_2({\ell}^n)$.
Then the rationality of $\boldsymbol{\sigma}$ implies that there exists an element $h_{\infty, p}$ in $G_2({\ell^n})$ with $\bar{g}_{\infty}^{h_{\infty, p}} = \bar{g}_{\infty}^p$.
It is not hard to see that relative to the Jordan canonical form of $\bar{g}_{\infty}$, $h_{\infty, p}$ is of the form:
\begin{equation}\label{con.eq.1}
\begin{pmatrix}
0 & 0 & 0 & 0 & 1 & 0 & 0 \\ 0 & 0 & 1 & 0 & 0 & 0 & 0 \\ 1 & 0 & 0 & 0 & 0 & 0 & 0 \\ 0 & 0 & 0 & 0 & 0 & 1 & 0 \\ 0 & 0 & 0 & 1 & 0 & 0 & 0 \\ 0 & 1 & 0 & 0 & 0 & 0 & 0 \\ 0 & 0 & 0 & 0 & 0 & 0 & 1
\end{pmatrix}^{\pm1}
\times {\rm ~an~element~of~} {\rm D}_7(\bar{\mathbb{F}}_\ell)
\end{equation}
where ${\rm D}_7$ denotes the diagonal torus of ${\rm GL}_7$.
Using a suitable Frobenius-action, we will construct in Rem.\ \ref{rem: h infty p II} below an element $h_{\infty,p} \in G_2({\ell})$ with the analogue property even for $p \neq 7$ but not necessary prime to the order of $G_2({\ell})$.
\end{rem}

\begin{rem}\label{rem: normalizer}
Let $K$ be any algebraically closed field. Then
\begin{equation*}
 N_{{\rm GL}_2(K)}(G_2(K)) = G_2(K) \times {\rm scalar~matrices~in~} K^\times:
\end{equation*}
Indeed, if $h$ is an element of the normalizer and $f(\cdot,\cdot,\cdot)$ denotes the Dickson alternating trilinear form then $G_2(K)$ still fixes the form $f(h(\cdot),h(\cdot),h(\cdot))$.
Hence by \cite{Aschbacher} Thm.\ 5 (2) this form is a scalar multiple of the original Dickson form.
Since $K$ is assumed to be algebraically closed, a scalar multiple of $h$ fixes the Dickson form.
By \cite{Aschbacher} (2.11) and (3.4) this means $h \in G_2(K) \times K^\times$ ($G_2$ is a simple group!).
The other inclusion is clear.
\end{rem}

\begin{rem}\label{rem: normalizer2}
Using $h_{\infty, p}$ from Rem.\ \ref{rem: h infty p} and the description of the normalizer of $G_2(\bar{\mathbb{F}}_{\ell})$ from Rem.\ \ref{rem: normalizer}, we get
\begin{equation*}
N_{{\rm GL}_7(\bar{\mathbb{F}}_{\ell})}(G_2({\ell^n})) = G_2({\ell^n}) \times {\rm scalar~matrices~in~} \bar{\mathbb{F}}_{\ell}^\times:
\end{equation*}
Obviously, $G_2({\ell}^n) \times \bar{\mathbb{F}}_{\ell}^\times$ normalizes $G_2(\ell^n)$, so we have to show the other inclusion.
For $g \in G_2(\ell^n)$ and $h \in N_{{\rm GL}_7(\bar{\mathbb{F}}_{\ell})}(G_2({\ell^n}))$ arbitrary we have $g^{\phi_{\ell}^n(h)} = \phi_{\ell}^n(g^h) = g^h$, i.e., $h^{-1}\phi_{\ell}^n(h)$ lies in the centralizer $C_{{\rm GL}_7(\bar{\mathbb{F}}_{\ell})}(G_2({\ell}))$, where $\phi_{\ell}$ is the component wise Frobenius map of $\bar{\mathbb{F}}_\ell / \mathbb{F}_\ell$.
Let $h^\prime$ be arbitrary in this centralizer.
We may assume that $\bar{g}_{\infty}$ is in its Jordan canonical form.
In particular $h^\prime$ is contained in ${\rm D}_7(\bar{\mathbb{F}}_\ell)$ (since this is the centralizer of $\bar{g}_\infty$).
And it follows from (\ref{con.eq.1}) that $h^\prime$ is of the form ${\rm diag}(a,\dots,a,b)$.
So e.g.\ from $\bar{g}_1^{h^\prime} = \bar{g}_1$ we can deduce that $a = b$.
In particular $h^{-1}\phi_{\ell}^n(h)$ lies in the center $\bar{\mathbb{F}}_\ell^\times$ of ${\rm GL}_7(\bar{\mathbb{F}}_{\ell})$.
Thus, there is an element $\lambda$ in this center s.t.\ $h^{-1}\phi_{\ell}^n(h)$ equals $\lambda\phi_\ell^n(\lambda)^{-1}$, i.e., the product $h \lambda$ is fixed by $\phi_\ell^n$.
It follows that $h \lambda$ is a $\phi_\ell^n$-fixed point of $N_{{\rm GL}_7(\bar{\mathbb{F}}_{\ell})}(G_2(\bar{\mathbb{F}}_{\ell})) = G_2(\bar{\mathbb{F}}_\ell) \times \bar{\mathbb{F}}_\ell^\times$, i.e., $h$ lies in $G_2({\ell}^n) \times \bar{\mathbb{F}}_{\ell}^\times$.
\end{rem}

\begin{lem}\label{lem: matzat rigid}
The tuple $\boldsymbol{\sigma}$ is rigid with respect to $G_2(\ell^n)$ in the sense of \cite{Matzat}.
\end{lem}

\proof
We just have to combine the linear rigidity of $\boldsymbol{\sigma}$ from Lem.\ \ref{lem: linear rigidity} with Rem.\ \ref{rem: normalizer2}.
\Endproof

Let us now prove the promised triviality of $n$:

\begin{cor}\label{cor: n=1}
Let ${\ell} \neq 2,3,7,11$ and $13$.
Then ${\rho}_{\bar{\mathcal{V}}}$ is already defined over $\mathbb{F}_\ell$ and the image of ${\rho}_{\bar{\mathcal{V}}}$ and $\rho_{\bar{\mathcal{V}}_{\mathbb{Q}}}$ is $G_2(\ell)$.
\end{cor}

\proof
The ${\rm GL}_7(\bar{\mathbb{F}}_{\ell})$-conjugacy classes of the $\bar{g}_i$'s are $\phi_{\ell}$-stable by Lem.\ \ref{lem: jordan forms remain untouched}.
Moreover, we still have the product relation
\begin{equation*}
 \phi_{\ell} (\bar{g}_0) \phi_{\ell} (\bar{g}_1) \phi_{\ell} (\bar{g}_{\infty}) = \textbf{1}.
\end{equation*}
Using Lem.\ \ref{lem: linear rigidity}, we get a single element $g \in {\rm GL}_7(\bar{\mathbb{F}}_{\ell})$ s.t.\ $\phi_{\ell}(\bar{g}_i)$ equals $\bar{g}_i^g$ for $i=0,1$ resp.\ $\infty$.
Since $\boldsymbol{\sigma}$ is a generating tuple of $G_2(\ell^n)$, $g$ lies in $N_{{\rm GL}_7(\bar{\mathbb{F}}_{\ell})}(G_2({\ell^n}))$, i.e., we may assume that $g$ lies in $G_2(\ell^n)$ by Rem.\ \ref{rem: normalizer2}.
Now $G_2$ is connected so there is an element $h \in G_2(\bar{\mathbb{F}}_\ell)$ with $g= h \phi_\ell(h)^{-1}$ by Lang-Steinberg (see e.g.\ \cite{Carter2} 1.17), i.e., $\phi_\ell$ fixes $\bar{g}_i^h$ while $\boldsymbol{\sigma}^h$ still is a generating tuple.
It follows that $\phi_\ell$ fixes $G_2(\ell^n)$ elementwise, i.e., $n=1$.
\Endproof

\section{Specializations of $\rho_{\bar{\mathcal{V}}_{\mathbb{Q}}}$}

From now on always assume ${\ell} \neq 2,3,7,11$ and $13$.
Denote the $\mathbb{Q}$-rational point of $\mathbb{A}_{\mathbb{Q}}^1 \setminus \{ 0,1 \}$ induced by $s \in \mathbb{Q} \setminus \{ 0,1 \}$ by $s$ as well.
Composition with the induced map on fundamental groups with $\rho_{\bar{\mathcal{V}}_{\mathbb{Q}}}$ gives us the specialization at $s$
\begin{equation*}
\xymatrix{
\rho_{\bar{\mathcal{V}}_{\mathbb{Q}}}^{(s)} : G_{\mathbb{Q}} \ar[r] & G_2(\ell)
}
\end{equation*}
(similar for $\mathcal{V}_k,\bar{\mathcal{V}}_k,\dots,$ where $k$ is a number field containing a primitive $7$th root of unity as above).
In this section we want to prove that the Galois representation $\rho_{\bar{\mathcal{V}}_{\mathbb{Q}}}^{(s)}$ is onto for suitable choices of $s$.

\subsection{Restriction of $\rho_{\bar{\mathcal{V}}_{\mathbb{Q}}}^{(s)}$ to inertia subgroups.}\label{sect: restriction to inertia subgroups}
First we have to modify $\bar{\mathcal{V}}_{R_N}, \bar{\mathcal{V}}_k$ s.t.\ the resulting sheaves are compatible with $\bar{\mathcal{V}}_{\mathbb{Q}}$.
Since $\pi_1^{{\rm \acute{e}t}}(\mathbb{A}_{\bar{\mathbb{Q}}}^1 \setminus \{ 0,1 \})$ is a normal subgroup of $\pi_1^{{\rm \acute{e}t}}(\mathbb{A}_{\mathbb{Q}(\zeta_7)}^1 \setminus \{ 0,1 \})$ and $\pi_1^{{\rm \acute{e}t}}(\mathbb{A}_{R_N}^1 \setminus \{ 0,1 \})$ is a quotient of the latter group, the image of $\rho_{\bar{\mathcal{V}}_{R_N}}$ is contained in $N_{{\rm GL}_7(\mathbb{F}(\zeta_7))}(G_2(\ell))$.
By Rem.\ \ref{rem: normalizer2} this normalizer is the product of $G_2(\ell)$ with the scalar matrices in $\mathbb{F}_{\ell}(\zeta_7)^\times$, i.e., we get a continuous character
$\varepsilon$ defined by the following composition of homomorphisms:
\begin{equation*}
\xymatrix{
\pi_1^{{\rm \acute{e}t}}(\mathbb{A}_{R_N}^1 \setminus \{ 0,1 \}) \ar[r]^-{\rho_{\bar{\mathcal{V}}_{R_N}}} \ar@/_1pc/[rrr]_{\varepsilon} & G_2({\ell}) \times 
\mathbb{F}_{\ell}(\zeta_7)^{\times} \ar[r]^-{{\rm pr}_2} & \mathbb{F}_{\ell}(\zeta_7)^{\times} \ar[r]^{(~)^{-1}} & \mathbb{F}_{\ell}(\zeta_7)^{\times}.
}
\end{equation*}
Note that the $7^{\rm th}$ $\otimes$-power of this character is just inverse to the character
given by the determinant of $\rho_{\bar{\mathcal{V}}_{R_N}}$: Indeed, $G_2(\ell)$ lies in ${\rm SL}_7(\ell)$, since the generators $\bar{g}_0,\bar{g}_1,\bar{g}_\infty$ do.
Define $\bar{\mathcal{V}}_{R_N}^{\varepsilon},\bar{\mathcal{V}}_{k}^{\varepsilon}, \dots$ as the twists of $\bar{\mathcal{V}}_{R_N},\bar{\mathcal{V}}_{k}, \dots$ by $\varepsilon$.
The restriction of $\varepsilon$ to the geometric fundamental group $\pi_1^{{\rm \acute{e}t}}(\mathbb{A}_{\bar{\mathbb{Q}}}^1 \setminus \{ 0,1 \})$ is trivial, so $\bar{\mathcal{V}}^{\varepsilon} = \bar{\mathcal{V}}$.
In particular these sheaves are still extensions of $\bar{\mathcal{V}}$, with
\begin{equation*}
 {\rm im}(\rho_{\bar{\mathcal{V}}_{k}^{\varepsilon}}) = G_2(\ell).
\end{equation*}
Now the restriction of $\bar{\mathcal{V}}_{\mathbb{Q}}$ to $\mathbb{A}_k^1 \setminus \{ 0,1 \}$ is a second extension having this property. It follows from the fact that the group $G_2(\ell)$ is centerfree that
 $\bar{\mathcal{V}}_k^{\varepsilon}$ is the restriction of $\bar{\mathcal{V}}_{\mathbb{Q}}$ to $\AA^1_k\setminus \{0,1\}.$

\begin{lem}\label{lem: specialization}
Let $p \neq 2,~7$ and ${\ell}$ be a prime number.
Let $i$ be $0$,$1$ or $\infty$.
Further, let $s \in \mathbb{Q} \setminus \{ 0,1 \}$ with $ \nu_p(s) >0 $ and  $2 \nmid \nu_p(s)$ (for $i = 0$), $ \nu_p(s-1) >0 $ and  ${\ell} \nmid \nu_p(s-1)$ (for $i = 1$) resp.\ $\nu_p(s) < 0$ and  $7 \nmid \nu_p(s)$ (for $i = \infty$).
Then the image of an inertia subgroup $ I_p \leq G_{\mathbb{Q}} $ under the specialization $\rho_{\bar{\mathcal{V}}_{\mathbb{Q}}}^{(s)}$ at $s$ is generated by a conjugate of $\bar{g}_i$.
\end{lem}

Before giving the proof, let us first recall a few facts on the local behaviour of \'{e}tale coverings of $\mathbb{A}^1\setminus \{ 0,1 \}$:

\begin{rem}\label{rem: local behaviour of coverings}
The absolute Galois groups resp.\ \'{e}tale fundamental group as well as the universal pro-\'{e}tale covering space of $\mathbb{C}_p (( z ))$ resp.\ $\Delta_p := {\rm Spec}(W_p[[z]][z^{-1}])$ are described explicitly in \cite{SGAI} XIII Cor.\ 5.3 applied to $\mathbb{C}_p[[z]]$ resp.\ $W_p[[z]]$ and the divisor given by $z$ (since the residue field of ${\rm div}(z)$ as a place of the respective function fields has characteristic $0$, every covering is tame).
Further, the explicit description of the universal pro-\'{e}tale covering space of $\Delta_p$ in loc.\ cit.\ implies that the map on fundamental groups given by evaluation by an $a\in W_p$
\begin{equation*}
 {\rm ev}_a:
 \xymatrix{
  W_p[[z]][z^{-1}] \ar[r] &
  \hat{\mathbb{Q}}_p^{{\rm nr}}
 }
\end{equation*}
factors over the tame quotient $\hat{\mathbb{Z}}^{(p^\prime)}(1)$ of the absolute Galois group of $\hat{\mathbb{Q}}_p^{{\rm nr}}$ (here $\hat{\mathbb{Z}}^{(p^\prime)}(1)$ means the prime-to-$p$ part).
We collect all these information in the following diagram:
\begin{equation*}
\xymatrix{
G_{\mathbb{C}_p((z))} \ar[r]^{\sim} \ar[d] & \hat{\mathbb{Z}}(1) \phantom{.} \ar[d]^{{\rm pr}}
\\
\pi_1^{{\rm \acute{e}t}}(\Delta_p) \ar[r]^{\sim} & \hat{\mathbb{Z}}^{(p^\prime)}(1) \phantom{.}
\\
G_{\hat{\mathbb{Q}}_p^{{\rm nr}}} \ar[r]^{{\rm pr}} \ar[u] & \hat{\mathbb{Z}}^{(p^\prime)}(1) . \ar[u]_{(-)^{\nu_p (a)}}
}
\end{equation*}
For the lower square, note that the pullback along ${\rm ev}_a$ of the covering $z \mapsto z^n$ of $\Delta_p$ is just the disjoint union of copies of the extension $K(\sqrt[d]p)/K$ for $d$ minimal s.t.\ $n$ divides $d\cdot \nu_p (a)$.
\end{rem}

We will apply Rem.\ \ref{rem: local behaviour of coverings} for $z = t$ (for $i = 0$), $z = t - 1$ (for $i = 1$) resp.\ $z = \frac{1}{t}$ (for $i = \infty$) and $a = z(s)$.\\
\phantom{leere Zeile}

\noindent \textbf{Proof of Lem.\ \ref{lem: specialization}:}
Choose an $N$ with $p \nmid N$.
The canonical inclusions of the resp.\ coordinate rings gives us the commutative diagram
\begin{equation*}
\xymatrix{
G_2({\ell}) &
\\
\pi_1^{{\rm \acute{e}t}}(\mathbb{A}_{R_N}^1 \setminus \{ 0,1 \}) \ar[u]_{\rho_{\bar{\mathcal{V}}_{R_N}^{\varepsilon}}}  & \pi_1^{{\rm \acute{e}t}}(\mathbb{A}_{\bar{\mathbb{Q}}}^1 \setminus \{ 0,1 \}) \ar[l] \ar[ul]_{\rho_{\bar{\mathcal{V}}}}
\\
\pi_1^{{\rm \acute{e}t}}(\Delta_p) \ar[u] & G_{\mathbb{C}_p((z))} \ar[l] \ar[u] \ar[r]^-{\sim} & G_{\bar{\mathbb{Q}}((z))}. \ar[ul]_{{\rm res}}
}
\end{equation*}
Write $\gamma_0,\gamma_1$ resp.\ $\gamma_\infty$ for the standard topological generators of $\pi_1^{{\rm \acute{e}t}}(\mathbb{A}_{\bar{\mathbb{Q}}}^1 \setminus \{ 0,1 \},\bar{s})$, i.e., $\gamma_0 \gamma_1 \gamma_\infty = \boldsymbol{1}$ and $\rho_{\bar{\mathcal{V}}}(\gamma_i)= \bar{g}_i$.
The induced map $G_{\bar{\mathbb{Q}}((z))} \rightarrow \pi_1^{{\rm \acute{e}t}}(\mathbb{A}_{\bar{\mathbb{Q}}}^1 \setminus \{ 0,1 \})$ identifies $G_{\bar{\mathbb{Q}}((z))}$ with the subgroup of $\pi_1^{{\rm \acute{e}t}}(\mathbb{A}_{\bar{\mathbb{Q}}}^1 \setminus \{ 0,1 \})$ topologically generated by $\gamma_i$ (look at the effect on convergent Laurent series of an \'{e}tale covering of $\mathbb{A}^1(\mathbb{C}) \setminus \{ 0,1 \}$ pulled back to open punctured discs around a missing point with radius converging to $0$).
It follows that the image of the resulting representation $G_{\bar{\mathbb{Q}}((z))} \rightarrow G_2(\ell)$ is the subgroup generated by $\bar{g}_i$.
The same is true for the representation $\pi_1^{{\rm \acute{e}t}}(\Delta_p) \rightarrow G_2(\ell)$ induced by $(\bar{\mathcal{V}}_{R_N})_{\mid \Delta_p}$:
Indeed, $G_{\mathbb{C}_p((z))} \rightarrow \pi_1^{{\rm \acute{e}t}}(\Delta_p)$ is an epimorphism by Rem.\ \ref{rem: local behaviour of coverings}.
\\
The canonical inclusions and evaluations $t \mapsto s$
give us the commutative diagram
\begin{equation*}
\xymatrix{
& & & G_2({\ell}) 
\\
\pi_1^{{\rm \acute{e}t}}(\mathbb{A}_{R_N}^1 \setminus \{ 0,1 \}, \bar{s}) \ar[rrru]^{\rho_{\bar{\mathcal{V}}_{R_N}^{\varepsilon}}} & & \pi_1^{{\rm \acute{e}t}}(\mathbb{A}_{\mathbb{Q}(\zeta_7)}^1 \setminus \{ 0,1 \}, \bar{s}) \ar[ll] \ar[ru]_{\rho_{\bar{\mathcal{V}}_{\mathbb{Q}(\zeta_7)}^{\varepsilon}}} \ar[r] & \pi_1^{{\rm \acute{e}t}}(\mathbb{A}_{\mathbb{Q}}^1 \setminus \{ 0,1 \}, \bar{s}) \ar[u]^{\rho_{\bar{\mathcal{V}}_{\mathbb{Q}}}}
\\
\pi_1^{{\rm \acute{e}t}}(\Delta_p, \bar{s}) \ar[u] & G_{\hat{\mathbb{Q}}_p^{{\rm nr}}} \ar[r]^-\sim \ar[ur] \ar[l] & I_p \ar@{^{(}->}[r] \ar[u] & G_{\mathbb{Q}}. \ar[u] \ar@/_3.8pc/[uu]_{\rho_{\bar{\mathcal{V}}_{\mathbb{Q}}}^{(s)}}
}
\end{equation*}
In particular, $(\rho_{\bar{\mathcal{V}}_{\mathbb{Q}}}^{(s)})_{\mid I_p}$ factors through the representation induced by $(\bar{\mathcal{V}}_{R_N})_{\mid \Delta_p}$.
By Rem.\ \ref{rem: local behaviour of coverings}, $G_{\hat{\mathbb{Q}}_p^{{\rm nr}}} \rightarrow \pi_1^{{\rm \acute{e}t}}(\Delta_p)$ is the $\nu_p (a)$ power map and $\nu_p (a)$ is coprime to ${\rm ord}(\bar{g}_i)$ by our assumptions.
It follows that the image of $(\rho_{\bar{\mathcal{V}}_{\mathbb{Q}}}^{(s)})_{\mid I_p}$ is still the subgroup generated by $\bar{g}_i$, which completes the proof.
\Endproof

Consider the maximal tame algebraic extension $\mathbb{Q}_p^{\rm tame} / \mathbb{Q}_p$.
The tame inertia group $I_p^{\rm tame}$ is given by the exact sequence
\begin{equation*}
\xymatrix{
\boldsymbol{1} \ar[r] & I_p^{\rm tame} \ar[r] & G_{\mathbb{Q}_p}^{\rm tame} \ar[r] & G_{\mathbb{F}_p} \ar[r] & \boldsymbol{1}.
}
\end{equation*}
Since $I_p^{\rm tame}$ is abelian, $G_{\mathbb{F}_p}$ acts on $I_p^{\rm tame}$.
It is not hard to see that the Frobenius acts as $p$ power map.
Now Lem.\ \ref{lem: specialization} provides us an element $h_{\infty,p}$ in the image of $\rho_{\bar{\mathcal{V}}_{\mathbb{Q}}}^{(s)}$ as in Rem.\ \ref{rem: h infty p}:

\begin{rem}\label{rem: h infty p II}
Let $p$ be as in the last lemma for $i = \infty$.
In particular, $p \neq 7$.
By Lem.\ \ref{lem: specialization}, the image of $(\rho_{\bar{\mathcal{V}}_{\mathbb{Q}}}^{(s)})_{\mid I_p}$ is generated by $\bar{g}_\infty$, i.e., is of order $7$.
Since $I_p^{\rm tame}$ is the quotient of $I_p$ by its unique $p$-Sylow group, $(\rho_{\bar{\mathcal{V}}_{\mathbb{Q}}}^{(s)})_{\mid I_p}$ factors through $I_p^{\rm tame}$.
In particular, a lift of the Frobenius to $G_{\mathbb{Q}_p}$ acts on the image of $(\rho_{\bar{\mathcal{V}}_{\mathbb{Q}}}^{(s)})_{\mid I_p}$ as the $p$ power map.
Thus for $p \equiv 3 ~{\rm or}~5~{\rm mod}~7$ we get an element $h_{\infty,p}$ in the image of $\rho_{\bar{\mathcal{V}}_{\mathbb{Q}}}^{(s)}$ as in Rem.\ \ref{rem: h infty p}.
\end{rem}

\subsection{The structure of $\bar{s}^*\bar{\mathcal{V}}_{\mathbb{Q}}$ as $\mathbb{F}_{\ell}[\bar{g}_1, \bar{g}_{\infty}, h_{\infty,p}]$-module.}
Choose a rational number $s$ satisfying the conditions of Lem.\ \ref{lem: specialization} for $i = 1$ and $\infty$ together with a suitable prime $p$ as in Rem.\ \ref{rem: h infty p II}.
Let $H$ be the subgroup of the image of $\rho_{\bar{\mathcal{V}}_{\mathbb{Q}}}^{(s)}$ generated by $\bar{g}_1, \bar{g}_{\infty}$ and $h_{\infty,p}$.
For a field extension $K / \mathbb{F}_{\ell}$ set $V(K) := \bar{s}^*\bar{\mathcal{V}}_{\mathbb{Q}} \otimes_{\mathbb{F}_{\ell}} K$.

\begin{lem}\label{im.lem.1.1}
The decomposition of $V(\bar{\mathbb{F}}_{\ell})$ into $V_1(\bar{\mathbb{F}}_{\ell}) \oplus V_2(\bar{\mathbb{F}}_{\ell})$, where the direct factors are $V_1(\bar{\mathbb{F}}_{\ell}) := {\rm Eig}(\bar{g}_{\infty},1)$ and $V_2(\bar{\mathbb{F}}_{\ell}) := \bigoplus_{i=1}^6 {\rm Eig}(\bar{g}_{\infty},\zeta_7^i)$, is defined over $\mathbb{F}_{{\ell}}$.
If $H$ acts reducible on $\bar{s}^*\bar{\mathcal{V}}_{\mathbb{Q}}$, the only possible non trivial $\mathbb{F}_{{\ell}}[H]$-submodules are $V_1(\mathbb{F}_{\ell})$ and $V_2(\mathbb{F}_{\ell})$.
Further, if they exist, they are irreducible submodules.
\end{lem}

\proof
The first statement is trivial.
For the second statement, consider the subgroup $H' \leq H$ generated by  $\bar{g}_{\infty}$ and $h_{\infty,p}$.
Let $v_i$ be a $\zeta_7^i$-eigenvector of $\bar{g}_{\infty}$.
Then (\ref{con.eq.1}) describes $h_{\infty,p}$ with resp.\ to the basis $\{ v_1, \dots, v_7\}$.
From this we see that $V(\bar{\mathbb{F}}_{\ell})=V_1(\bar{\mathbb{F}}_{\ell}) \oplus V_2(\bar{\mathbb{F}}_{\ell})$ is a decomposition of $V(\bar{\mathbb{F}}_{\ell})$ as a $\bar{\mathbb{F}}_{\ell}[H']$-module.
We claim that it is even an irreducible decomposition of $V(\bar{\mathbb{F}}_{\ell})$ as a $\bar{\mathbb{F}}_{\ell}[H']$-module, i.e., it remains to prove the $H'$-irreducibility of $V_2(\bar{\mathbb{F}}_{\ell})$:
For a non trivial $v \in V_2(\bar{\mathbb{F}}_{\ell})$ define $M_v$ as the set of all $i$ s.t.\ the projection of $v$ to $\bar{\mathbb{F}}_{\ell}.v_i$ is non trivial and $m_v$ as $\vert M_v \vert$.
Let $P \in \bar{\mathbb{F}}_{\ell}[X]$ be non trivial of minimal degree with $P(\bar{g}_{\infty})v = 0$.
It follows that $P(\zeta_7^i) = 0$ for $i \in M_v$.
Now $\prod_{i \in M_v} (X - \zeta_7^i)$ gives another such non trivial relation, i.e., ${\rm deg}(P) = m_v$.
The sum $\sum_{i \in \mathbb{Z}} \bar{\mathbb{F}}_{\ell}.\bar{g}_{\infty}^i (v)$ is contained in the $m_v$-dimensional space $\bigoplus_{i\in M_v} \bar{\mathbb{F}}_{\ell}.v_i$.
As a result these spaces coincide.
In particular, there exists an $ i \leq 6$ with $v_i$ in the $\bar{\mathbb{F}}_{\ell}$-span of $H'v$.
Using (\ref{con.eq.1}) we conclude that this span is already $V_2(\bar{\mathbb{F}}_{\ell})$.
\\
Let now $U \leq V(\mathbb{F}_{\ell})$ be a non trivial $\mathbb{F}_{\ell}[H]$- hence also $\mathbb{F}_{\ell}[H']$-submodule.
Using the $H'$-irreducibility of $V_1(\mathbb{F}_{\ell})$ and $V_2(\mathbb{F}_{\ell})$ we get that $U$ equals either $V_1(\mathbb{F}_{\ell})$ or $V_2(\mathbb{F}_{\ell})$ (compute the possible dimensions of their intersections with $U$).
\Endproof

It follows that $\bar{s}^*\bar{\mathcal{V}}_{\mathbb{Q}} =V_1(\mathbb{F}_{\ell}) \oplus V_2(\mathbb{F}_{\ell})$ is the only possible non trivial decomposition as an $\mathbb{F}_{\ell}[H]$-module.
Suppose that this is indeed a decomposition.
In particular both $V_1(\mathbb{F}_{\ell})$ and $V_2(\mathbb{F}_{\ell})$ are $\bar{g}_1$-invariant, which contradicts that the Jordan canonical form of $\bar{g}_1$ contains no Jordan block of length 1.
Thus we get:

\begin{cor}\label{im.cor.1.2}
$\bar{s}^*\bar{\mathcal{V}}_{\mathbb{Q}}$ is indecomposable as an $\mathbb{F}_{{\ell}}[H]$-module.
\end{cor}

\subsection{Computation of the image of $\rho_{\bar{\mathcal{V}}_{\mathbb{Q}}}^{(s)}$.}\label{image}
First we recall the following fact about our representation of $G_2({\ell})$ which can be found for example in \cite{Malle}:
There is a basis $\{ v_1, \dots, v_7 \}$ of $\bar{s}^*\bar{\mathcal{V}}_{\mathbb{Q}}$ s.t.\
every maximal parabolic subgroup of $G_2({\ell})$ is given either by $P_{\alpha}:=\langle \mathcal{T}, x_{\alpha}, x_{\alpha}^{\rm tr}, x_{\beta} \rangle$ or $P_{\beta}:=\langle \mathcal{T}, x_{\alpha}, x_{\beta}, x_{\beta}^{\rm tr} \rangle$, where $x_{\alpha}$ respectively $x_{\beta}$ are generators of the root subgroup of the simple root $\alpha$ respectively $\beta$ given as
\begin{equation*}
x_{\alpha}:=
\begin{pmatrix}
1 & 1 & & & & & \\ & 1 & & & & & \\ & & 1 & 1 & -1 & & \\ & & & 1 & -2 & & \\ & & & & 1 & & \\& & & & & 1 & -1 \\ & & & & & & 1
\end{pmatrix},~x_{\beta}:=
\begin{pmatrix}
1 & & & & & & \\ & 1 & 1 & & & & \\ & & 1 & & & & \\ & & & 1 & & & \\ & & & & 1 & -1 & \\& & & & & 1 & \\ & & & & & & 1
\end{pmatrix}
\end{equation*}
and $\mathcal{T} := \lbrace {\rm diag}(t_1,t_2,t_1t_2^{-1},1,t_1^{-1}t_2,t_2^{-1},t_1^{-1}) \mid t_1,t_2 \in \mathbb{F}_{{\ell}}^{\times} \rbrace$ is a maximal split torus of $G_2({\ell})$ in ${\rm GL}_7(\mathbb{F}_{{\ell}})$.

We will use the above list of types of maximal subgroups of $G_2(\ell)$ to prove the following

\begin{thm}\label{thm: hilbert irreducibility}
Let ${\ell} \neq 2,3,7,11,13$ be a prime number.
Let $s \in \mathbb{Q} \setminus \lbrace 0,1 \rbrace$ be a rational number and suppose there are primes $p, q \neq 2, 7, {\ell}$ with $p \equiv 3 ~{\rm or}~ 5 ~{\rm mod}~ 7$, $\nu_p(s) < 0 < \nu_q(s - 1)$, $7 \nmid \nu_p(s)$ and ${\ell} \nmid \nu_q(s - 1)$.
Then $\rho_{\bar{\mathcal{V}}_{\mathbb{Q}}}^{(s)}: G_{\mathbb{Q}} \rightarrow G_2({\ell})$ is surjective.
\end{thm}

Before we give the proof, let us make a few comments about the theorem:

\begin{rem}
Note that Thm.\ \ref{thm: hilbert irreducibility} is in particular an explicit version of Hilbert's irreducibility theorem for $G_2({\ell})$.
One could get a similar result from S.\ Beckmann's work in \cite{Beckmann}: to apply Cor.\ 1.3 in loc.\ cit.\ to $G_2(\ell)$ one needs a generating set of certain $\lbrace \sigma_i \mid i \in T \rbrace$ in the image of certain inertia groups with the property that for all $\tau_i$ in the $G_2({\ell})$-conjugacy class of $\sigma_i$, $\lbrace \tau_i \mid i \in T \rbrace$ is still a generating set.
But restriction of our specialization to the inertia groups $I_p$ and $I_q$ alone only yields two elements, which not necessary generates $G_2({\ell})$.
Thus her result does not apply in our case.
Further, we have to exclude less primes: $p$ and $q$ might as well divide the order of $G_2({\ell})$, hence might lie inside the exceptional set $S_{{\rm bad}}$ in \cite{Beckmann}.
\end{rem}

\begin{rem}\label{rem: infinitely many s}
There are infinitely many $s\in \mathbb{Q}\setminus \{0,1\}$ satisfying the assumptions of Thm.\ \ref{thm: hilbert irreducibility}:
First, note that there are infinitely many primes $p, q$  s.t.\ $p, q \neq 2, 7, {\ell}$ and $p \equiv 3 ~{\rm or}~ 5 ~{\rm mod}~ 7$ (for the last condition use Dirichlet's prime number theorem).
Now for each $a_p, a_q \in \mathbb{Q}$ s.t. $ \nu_p(a_p) < 0$ and $7 \nmid \nu_p(a_p)$ resp.\ $\nu_q(a_q - 1) > 0$ and $\ell \nmid \nu_q(a_q - 1)$ we find an $s\in \mathbb{Q}$ arbitrary close to $a_p$ resp.\ $a_q$ with respect to the $p$- resp.\ $q$-adic valuation by the approximation theorem.
Since there are infinitely many such pairs of $a_p$ and $a_q$ with different $p$- resp.\ $q$-adic valuations, we get infinitely many $s$.  
\end{rem}

\begin{rem}\label{rem: assumptions}
The following equivalent formulation of the assumptions of Thm.\ \ref{thm: hilbert irreducibility} is better suited for Sect.\ \ref{sect: motives}:
Let $s \in \mathbb{Q} \setminus \lbrace 0,1 \rbrace$ be a rational number and suppose there are primes $p, q \neq 2, 7$ with $p \equiv 3 ~{\rm or}~ 5 ~{\rm mod}~ 7$, $\nu_p(s) < 0 < \nu_q(s - 1)$ and $7 \nmid \nu_p(s)$.
Then $\rho_{\bar{\mathcal{V}}_{\mathbb{Q}}}^{(s)}$ is surjective for all primes $\ell \neq 2,3,7,11,13,p,q$ s.t.\ ${\ell} \nmid \nu_q(s - 1)$.
In particular, for a fixed $s$ satisfying the first set of conditions, the theorem is true for all but finitely many primes $\ell$. 
\end{rem}

\noindent \textbf{Proof of Thm.\ \ref{thm: hilbert irreducibility}:}
As above, let $H$ be the subgroup of the image of $\rho_{\bar{\mathcal{V}}_{\mathbb{Q}}}^{(s)}$ generated by $\bar{g}_1, \bar{g}_{\infty}$ and $h_{\infty,p}$.
It suffices to prove that $H$ is not contained inside any of the maximal subgroups $G_2({\ell})$-conjugate to one of the subgroups in the list of Rem.\ \ref{rem: maximal subgroups}.
\\
By the same arguments as in the proof of Prop.\ \ref{prop: image G2}, $H$ is not contained in any maximal subgroup of $G_2({\ell})$ (conjugate to a subgroup) in the above list except the first three ones.
\\
\textbf{(i)}
Suppose $H \leq M$ with $M$ a maximal parabolic subgroup of $G_2({\ell})$.
If $M=P_{\alpha}$, then $v_1$ and $v_2$ span a 2-dimensional $\mathbb{F}_{{\ell}}[P_{\alpha}]$- hence $\mathbb{F}_{{\ell}}[H]$-submodule of $\bar{s}^*\bar{\mathcal{V}}_{\mathbb{Q}}$.
If $M=P_{\beta}$, then $v_1,v_2$ and $v_3$ span a 3-dimensional $\mathbb{F}_{{\ell}}[P_{\beta}]$- hence $\mathbb{F}_{{\ell}}[H]$-submodule of $\bar{s}^*\bar{\mathcal{V}}_{\mathbb{Q}}$, which both contradicts Lemma \ref{im.lem.1.1}.
\\
\textbf{(ii)}
Suppose $H \leq C_{G_2({\ell})}(\iota)$ for an involution $\iota$ as in the above list.
Now $\iota$ is semisimple $\neq -\textbf{1}$ (otherwise $C_{G_2({\ell})}(\iota)$ would be $G_2({\ell})$).
Thus $\iota$ has Jordan canonical form
\begin{equation*}
{\rm JCF}(\iota) =
\begin{pmatrix}
{\rm id}_{r \times r} & 0 \\ 0 & -{\rm id}_{s \times s}
\end{pmatrix}
\end{equation*}
for suitable $r,s > 0$ and $r+s = 7$ (one has $r=3,s=4,$ cf.~\cite{Dettweiler1}, Table~1). 
It follows that $C_{{\rm GL}(\bar{s}^*\bar{\mathcal{V}}_{\mathbb{Q}})}(\iota)$ is isomorphic to ${\rm GL}_r ({\ell}) \oplus {\rm GL}_s ({\ell})$.
In particular even $C_{{\rm GL}(\bar{s}^*\bar{\mathcal{V}}_{\mathbb{Q}})}(\iota)$ acts decomposably on $\bar{s}^*\bar{\mathcal{V}}_{\mathbb{Q}}$, which contradicts Cor.\ \ref{im.cor.1.2}.
\\
\textbf{(iii)}
Suppose $H \leq M$ with $M$ $G_2({\ell})$-conjugate to $K_{\varepsilon} = L_{\varepsilon} \rtimes \mathbb{Z} / 2$.
Without loss of generality we may assume that $M$ equals $K_{\varepsilon}$.
According to the proof of \cite{Kleidman} Prop.\ 2.2 we may choose $L_{\varepsilon}$ in our representation to be generated by three long root subgroups 
of $G_2({\ell})$, all of which act trivially on a 5-dimensional subspace $W_i\, (i=1,2,3)$ of $\bar{s}^*\bar{\mathcal{V}}_{\mathbb{Q}}$ (resp.).
In particular $L_{\varepsilon}$ acts trivially on the non trivial subspace $W_1 \cap W_2 \cap W_3$.
Say $v$ is a non trivial element of this space and let $\sigma$ be the generator of $\mathbb{Z} / 2$ in $K_{\varepsilon}$.
Since $\sigma$ normalizes $L_{\varepsilon}$, $L_{\varepsilon}$ acts trivially on $\sigma v$ as well.
Thus the space $W$ spanned by $v$ and $\sigma v$ is a 1- or 2-dimensional $\mathbb{F}_{{\ell}}[K_{\varepsilon}]$-submodule of $\bar{s}^*\bar{\mathcal{V}}_{\mathbb{Q}}$.
Further, \cite{Kleidman} Prop.\ 2.2 (iii) tells us that $K_{\varepsilon}$ stabilizes a 6-dimensional subspace $U$ of $\bar{s}^*\bar{\mathcal{V}}_{\mathbb{Q}}$.
From Lem.\ \ref{im.lem.1.1} we get that $U = V_2(\ell)$ is an irreducible $\mathbb{F}_{{\ell}}[H]$-submodule.
In particular $U$ and $W$ intersect trivially.
It follows that $\bar{s}^*\bar{\mathcal{V}}_{\mathbb{Q}} = U \oplus W$ is a decomposition as an $\mathbb{F}_{\ell}[H]$-module, which again contradicts Cor.\ \ref{im.cor.1.2}, proving the theorem.\Endproof

\section{Motives with motivic Galois group $G_2$}\label{sect: motives}

Recall from Sect.\ \ref{basicob} that our lisse sheaf $\mathcal{V}$ on $\mathbb{A}_{\bar{\mathbb{Q}}}^1 \setminus \{0,1\}$ was the restriction of a lisse sheaf on $\mathbb{A}_{R_N}^1 \setminus \{0,1\}$ for $R_N= \mathbb{Z}[\zeta_7,\frac{1}{14\cdot l\cdot N}]$.
As such, we get canonical extensions $\mathcal{V}_k$ to $\mathbb{A}_{k}^1 \setminus \{0,1\}$ for any number field containing $\zeta_7$.
Out of $\mathcal{V}_k$ we will construct a family of motives with motivic Galois group (with resp.\ to the $\ell$-adic realization) $G_2(\bar{\mathbb{Q}}_{\ell})$ over a suitable quadratic field extension $k / \mathbb{Q}(\zeta_7)$.
But first, let $k / \mathbb{Q}(\zeta_7)$ be any algebraic field extension s.t.\ the Galois action given by the restriction of $\rho_{\bar{\mathcal{V}}_{\mathbb{Q}}}^{(s)}$ is non trivial (e.g.\ $k = \mathbb{Q}(\zeta_7)$ itself).

\subsection{Interpretation of $\bar{s}^*\mathcal{V}_k$ as $\ell$-adic realization.}
From the construction in section \ref{basicob} the way is clear: we have to interpret suitable stalks of $\mathcal{V}_k$ in terms of $\ell$-adic cohomology of suitable smooth projective varieties.

First, observe that the generic fibre ${\rm Hyp}_{k}$ in $\mathbb{G}_{{\rm m},{k}} \otimes_{k} (\mathbb{A}_{k}^7 \setminus v(\Delta))$ of ${\rm Hyp}_N$ (see Sect.\ \ref{basicob}) is already defined over $\mathbb{Q}$.
Denote the corresponding geometrically connected $\mathbb{Q}$-variety by ${\rm Hyp}_{\mathbb{Q}}$.
Although it is not \'{e}tale Galois over $\mathbb{A}_{\mathbb{Q}}^7 \setminus v(\Delta)$, it is still an \'{e}tale covering space.

\begin{prop}\label{prop: motivic interpretation in families}
Let $S$ be $\mathbb{A}^1 \setminus \{ 0,1 \}$.
There exists a smooth projective compactification $X_{\mathbb{Q}}/S$ of ${\rm Hyp}_{\mathbb{Q}}/S$ whose complement is a simple normal crossing divisor $D_{\mathbb{Q}} = \sum_i D_{i,\mathbb{Q}}$ over $S$ and such that  the base change
\begin{equation*}
\xymatrix{
{\rm Hyp}_{k} \ar@{^{(}->}[r]^j \ar[dr]_{\pi_K} & X \ar[d]^{\pi_{X}} & D \ar@{_{(}->}[l]_i \ar[dl]^{\pi_{D}}
\\
& \mathbb{A}_{k}^1 \setminus \{0,1\}
}
\end{equation*}
to ${k}$ is $\boldsymbol{\mu}_{14}$-equivariant.
Further, the canonical map $\iota:\coprod_i D_i \rightarrow X$ with induced projection $\pi_{\coprod D}$ induces an isomorphism
\begin{equation*}
 \mathcal{V}_{k} \cong \Pi[{\rm ker}(\mathds{R}^6\pi_{X,*}\bar{\mathbb{Q}}_\ell \rightarrow \mathds{R}^6\pi_{\coprod D,*}\bar{\mathbb{Q}}_\ell)],
\end{equation*}
where $\Pi$ is the idempotent endomorphism given in the group ring $\bar{\mathbb{Q}}_{\ell}[\boldsymbol{\mu}_{14}]$ as $\frac{1}{14} \sum_{i = 0}^{13} \zeta_{14}^{-i} \omega^i$ for $\omega$ the generator of $\boldsymbol{\mu}_{14}$ satisfying $\chi(\omega) = \zeta_{14}$.
\end{prop}

\proof
Let $\mathbb{P}_S^6$ the relative projective space with homogeneous coordinates $X_0$, \dots, $X_6$ and let $X_7$ be the coordinate of $S = \mathbb{A}^1 \setminus \{ 0,1 \}$.
The Zariski closure $X_{\mathbb{Q}}^{(0)}$ of ${\rm Hyp}_{\mathbb{Q}}$ in $\mathbb{P}^1 \times \mathbb{P}_S^6$ is a ramified covering of $\mathbb{P}_S^6$ whose (reduced) ramification locus is contained inside the hyperplane arrangement given by $\Delta$ and the hyperplane at infinity $X_0=0$.
Note that the multiplicities of the ramification divisor are constant over $S$.
As in the proof for \cite{Dettweiler1} Cor.\ 2.4.2 we use the standard resolution of a linear hyperplane arrangement given in \cite{EsnaultSchechtmanViehweg} Section~2 to get a projective birational map $P_S^{(1)} \rightarrow \mathbb{P}_S^6$ s.t.\ the pullback of the above hyperplane arrangement is a simple normal crossing divisor over $S$ and its strict transform is non singular.
\\
Let $X_{\mathbb{Q}}^{(1)}$ be the pullback of $X_{\mathbb{Q}}^{(0)}$ to $P_S^{(1)}$.
It is a ramified covering of $P_S^{(1)}$ whose (reduced) ramification locus is a strict normal crossing divisor over $S$ and whose ramification divisor still has constant multiplicities over $S$.
Thus, the singularities of this covering space are locally given by equations $Y^{14} = Z_0^{a_0} \cdots Z_6^{a_6}$ (see the equation for ${\rm Hyp}_N$ in Sect.\ \ref{basicob} for the $Y$-exponent) with constant multiplicities $a_i$ over $S$ and $Z_i$ of the form $X_i + f_i(X_7)$.
It follows that the fibres $s^*X_{\mathbb{Q}}^{(1)}$ locally are toric varieties whose fans $\Sigma = \Sigma_s$ are independent from $s\in S$.
Choose a refinement $\Sigma^\prime$ of $\Sigma$ as in \cite{CoxLittleSchenck} Thm.\ 11.2.2 corresponding to simple normal crossing resolutions of $X_{\Sigma_s}$.
This simple normal crossing resolutions of $X_{\Sigma_s}$ is purely given by the combinatorics of the common refinement $\Sigma^\prime$ of $\Sigma_s$. 
Therefore, the corresponding series of blowups are compatible and define a simple normal crossing resolution $X_{\mathbb{Q}}^{(2)} \rightarrow X_{\mathbb{Q}}^{(1)}$ of a given singularity in the whole family.
Iterating these locally toric resolutions for the remaining singularities in $X_{\mathbb{Q}}^{(2)}$ (still of the above local form), we get a series of analogue resolutions $X_{\mathbb{Q}}^{(m+1)} \rightarrow X_{\mathbb{Q}}^{(m)}$ with $X_{\mathbb{Q}} := X_{\mathbb{Q}}^{(n)}/S$ smooth for $n \gg 0$.
\\
The base change $X_k^{(1)}$ of $X_{\mathbb{Q}}^{(1)}$ to $k$ is a ramified $\boldsymbol{\mu}_{14}$-covering of $P_k^{(1)}$.
From the above local description of $X_k^{(m)}$ as a toric variety we get that $\boldsymbol{\mu}_{14}$ acts as a subgroup of the corresponding torus.
Toric resolutions are equivariant under the torus-action, so $X:= X_{\mathbb{Q}}\otimes k$ is our desired equivariant simple normal crossing compactification.
\\
The remaining claim follows from \cite{WeilII} Th\'eor\`eme~1 (Weil II) applied to the excision sequence of a suitable integral model of this compactification. For the detailed arguments, see the proof of \cite{Dettweiler1} Cor.\ 2.4.2.
\Endproof

Let $L/k$ be any field extension and $s\in S(L)$ an $L$-point.
As in \cite{Dettweiler1} Sect.\ 3.3, we get a motive for motivated cycles (see \cite{Andre96})
\begin{equation*}
 N_s := p_s\mathfrak{h}(X_s)
\end{equation*}
over $L$ with coefficients in $E := \mathbb{Q}(\zeta_7)$, where $p_s$ is the composition of the commuting motivated correspondences $\Pi_s$, ${\rm pr}_{{\rm ker}(\iota_s)}$ and the sixth K\"{u}nneth projector $\pi_{X_s}^6$.
Further, for $L=\mathbb{C}$ these motives $N_s$ form a family of motives over $S(\mathbb{C})$ in the sense of \cite{Andre96} Sect.\ 5.2.
We define the motive $M_s$ as the third Tate-twist $N_s(3)$.
Finally, let $G_{N_s}^{\rm mot}, G_{M_s}^{\rm mot}, \dots$ be the ($E$-linear) motivic Galois groups of $N_s, M_s. \dots$ as motives for motivated cycles.

\begin{cor}\label{Cor: motivic interpretation}
The ${\ell}$-adic realization $H_{\ell}^\bullet(M_s)$ of $M_s$ is isomorphic to $\bar{s}^*\mathcal{V}_k$ as a $G_k$-module.
\end{cor}

Here, by the $\ell$-adic realization $H_\ell^\bullet(-)$, we mean $H_{\rm \acute{e}t}^\bullet(-,\bar{\mathbb{Q}}_\ell)$.

\begin{rem}\label{rem: numerical motive}
If $X$ satisfies the K\"{u}nneth conjecture, then $\pi_X^6$ is algebraic and $p_s$ an idempotent in the category $\underline{\rm NM}(k)_{E}$ of numerical motives.
Thus, $M_s$ is a well defined numerical motive in this case.
Under the assumption of all the standard conjectures on algebraic cycles, the category  $\underline{\rm NM}(k)_{E}$ agrees with the Tannakian category of motives for motivated cycles and $G_{M_s}^{\rm mot}$ is the motivic Galois group in the sense of Grothendieck.
\end{rem}

\begin{rem}\label{rem: motive defined over the rationals}
If $s$ is a $\mathbb{Q}$-rational point, then ${\rm pr}_{{\rm ker}(\iota_s)}$ and $\pi_{X_s}^6$ are even rational motivated cycles, so $M_s^0 := {\rm pr}_{{\rm ker}(\iota_s)} \pi_{X_s}^6 \mathfrak{h}(X_s)(3)$ is the base change of a rational motive over $\mathbb{Q}$ to $\mathbb{Q}(\zeta_7)$ and $M_s = \Pi_s M_s^0$ is a motive over $\mathbb{Q}(\zeta_7).$
It follows that the ${\ell}$-adic realization of $M_s$ is just the $\zeta_7$-eigenspace of $\omega$ acting on the $\ell$-adic realization $H_{\ell}^\bullet(M_s^0)$ of the rational motive $M_s^0$.
\end{rem}

\subsection{The Galois-action on $H_\ell^\bullet(M_s)$.}
In this section we will prove that for a suitable quadratic field extension $k/\mathbb{Q}(\zeta_7)$, suitable choices of $s \in \mathbb{Q}\setminus \{ 0,1 \}$ and a prime $\ell \gg 0$ the absolute Galois group $G_k$ acts on the $\ell$-adic realization $H_\ell^\bullet(M_s)$ as a Zariski dense subgroup of $G_2$.
\\
So let $s \in \mathbb{Q} \setminus \lbrace 0,1 \rbrace$ be a rational number and suppose there are prime numbers $p, q \neq 2, 7$ with $p \equiv 3 ~{\rm or}~ 5 ~{\rm mod}~ 7$, $\nu_p(s) < 0 < \nu_q(s - 1)$ and $7 \nmid \nu_p(s)$.
Fix a prime $\ell \neq 2,3,7,11,13,p,q$ s.t.\ ${\ell} \nmid \nu_q(s - 1)$.

Since $\rho_{\mathcal{V}} $ 
has Zariski dense image in $G_2(\bar{\mathbb{Q}}_\ell)$, the image of 
$\rho_{{\mathcal{V}}_{R_N}(m)}$ (with $\mathcal{V}_{R_N}(m)$ denoting  
the $m$th Tate twist of $\mathcal{V}_{R_N}$ as usual) 
normalizes $G_2(\bar{\mathbb{Q}}_\ell)$ (see the first remarks in Sect.\ \ref{sect: restriction to inertia subgroups}).
Again, for $m=0$, the resulting character (use Rem.\ \ref{rem: normalizer})
\begin{equation*}
\varepsilon:
\xymatrix{
\pi_1^{{\rm \acute{e}t}}(\mathbb{A}_{R_N}^1 \setminus \{ 0,1 \}) \ar[r]^-{\rho_{{\mathcal{V}}_{R_N}}} \ar@/_1pc/[rrr] & G_2(\bar{\mathbb{Q}}_{\ell}) \times \bar{\mathbb{Q}}_{\ell}^{\times} \ar[r]^-{{\rm pr}_{\bar{\mathbb{Q}}_{\ell}^{\times}}} & \bar{\mathbb{Q}}_{\ell}^{\times} \ar[r]^{(~)^{-1}} & \bar{\mathbb{Q}}_{\ell}^{\times}
}
\end{equation*}
is just inverse to a character whose $7^{\rm th}$ $\otimes$-power is the determinant of $\rho_{\mathcal{V}_{R_N}}$.
Now, the restriction of $\varepsilon$ to $\mathbb{A}_{\bar{\mathbb{Q}}}^1 \setminus \{ 0,1 \}$ is trivial, so the twists $\mathcal{V}_k^{\varepsilon}, \mathcal{V}_{R_N}^{\varepsilon}, \dots$ of $\mathcal{V}_k, \mathcal{V}_{R_N}, \dots$ by the restricted character $\varepsilon$ are extensions of $\mathcal{V}$ with monodromy inside $G_2(\bar{\mathbb{Q}}_\ell)$.
\\
Note that the mod-$\ell$-reduction of $\varepsilon$ (the character is still defined over $\mathbb{Q}_{\ell}(\zeta_7)$) is just the character $\varepsilon$ defined as in Sect.\ \ref{sect: restriction to inertia subgroups} (justifying the same notation).
In particular, we get back $\bar{\mathcal{V}}_k^{\varepsilon}, \bar{\mathcal{V}}_{R_N}^{\varepsilon}, \dots$ as the mod $\ell$ reduction of $\mathcal{V}_k^{\varepsilon}, \mathcal{V}_{R_N}^{\varepsilon}, \dots$ in analogy to the mod $\ell$ reductions $\bar{\mathcal{V}}_k, \bar{\mathcal{V}}_{R_N}, \dots$ of $\mathcal{V}_k, \mathcal{V}_{R_N}, \dots$ at the end of Sect.~\ref{basicob}.

\begin{lem}\label{lem: irreducibility of the tate realization}
Let $k$ be an extension of  $\mathbb{Q}(\zeta_7)$ of degree $\leq 2.$
Then $\bar{s}^*\mathcal{V}_k(m) \cong H_{\ell}^\bullet(N_s(m))$ is an irreducible $G_k$-module.
\end{lem}

\proof
The $G_k$-module $\bar{s}^*\mathcal{V}_k(m)$ is the twist of $\bar{s}^*\mathcal{V}_k^{\varepsilon}$ by the character $\varepsilon^{-1}(m)$.
Thus, it remains to show that $\bar{s}^*\mathcal{V}_k^{\varepsilon}$ is irreducible:
Since $\rho_{\bar{\mathcal{V}}_k^{\varepsilon}}$ is the restriction of the absolutely irreducible $G_{\mathbb{Q}}$-representation $\rho_{\bar{\mathcal{V}}_{\mathbb{Q}}}$ (see Thm.\ \ref{thm: hilbert irreducibility}) with image the simple group $G_2(\ell)$, $G_k$ acts on $\bar{s}^*\bar{\mathcal{V}}_k^{\varepsilon}$ absolutely irreducibly, too.
It follows that the $G_k$-module $\bar{s}^*\mathcal{V}_k^{\varepsilon}$ is irreducible, as well.
\Endproof

\begin{lem}\label{lem: choice of the quadratic extension}
There is an extension $k/\mathbb{Q}(\zeta_7)$ of degree $\leq 2$ s.t.\ the restriction of $\varepsilon(-3)$ to $\mathbb{A}_k^1\setminus \{0,1\}$ is trivial.
In particular, $\mathcal{V}_k^{\varepsilon} = \mathcal{V}_k(3)$.
\end{lem}

\proof
Start with $k = \mathbb{Q}(\zeta_7)$.
The dimension of $X$ is $6$.
Poincar\'{e}-Duality gives us a $G_{\mathbb{Q}(\zeta_7)}$-equivariant non degenerate symmetric bilinear form $H_{\ell}^6(X)(3)^{\otimes 2} \rightarrow \bar{\mathbb{Q}}_{\ell}$.
Further, $M_s$ is defined as $p_s \mathfrak{h}(X)(3)$ and $H_{\ell}^\bullet(p_s)$ is $G_{\mathbb{Q}(\zeta_7)}$-equivariant.
It follows that the kernel of the canonical map $H_{\ell}^\bullet(M_s) \rightarrow H_{\ell}^\bullet(M_s)^{\vee}$ induced by the above bilinear form is $G_{\mathbb{Q}(\zeta_7)}$-invariant.
By Lem.\ \ref{lem: irreducibility of the tate realization}, $G_{\mathbb{Q}(\zeta_7)}$ acts irreducibly on $H_{\ell}^\bullet(M_s)$, so the restriction of the above bilinear form is a $G_{\mathbb{Q}(\zeta_7)}$-equivariant non degenerate symmetric bilinear form $b$ on $H_{\ell}^\bullet(M_s)$ resp.\ $\bar{s}^*\mathcal{V}_k(3)$.
We conclude that the image of $\rho_{\mathcal{V}_k(3)}^{(s)}$ is even contained inside
\begin{equation*}
 N_{{\rm O}(\bar{s}^*\mathcal{V}_k(3),b)}(G_2(\bar{\mathbb{Q}}_{\ell})) = G_2(\bar{\mathbb{Q}}_{\ell}) \times \langle \pm 1 \rangle.
\end{equation*}
In particular, $\varepsilon(-3)^{\otimes 2}$ is trivial and we can choose an extension $k/\mathbb{Q}(\zeta_7)$ of degree $\leq 2$ in Lem.\ \ref{lem: irreducibility of the tate realization} with the additional property $\mathcal{V}_k^{\varepsilon} = \mathcal{V}_k(3)$.
\Endproof

For the rest of the paper, we fix $k/\mathbb{Q}(\zeta_7)$ as in Lem \ref{lem: choice of the quadratic extension}.

\begin{thm}\label{thm: tate realization}
Let $s\in \mathbb{Q}\setminus \{ 0,1 \}$ be a rational number and suppose there are prime numbers $p,q \neq 2, 7$ s.t.\ $p \equiv 3 ~{\rm or}~ 5 ~{\rm mod}~ 7$, $\nu_p(s) < 0 < \nu_q(s - 1)$ and $7 \nmid \nu_p(s)$.
Let $k/\mathbb{Q}(\zeta_7)$ be the extension of degree $\leq 2$ given in Lem \ref{lem: choice of the quadratic extension}.
We identify $s$ with the induced $k$-point of $\mathbb{A}^1\setminus \{ 0,1 \}$.
Then for any prime ${\ell} \neq 2,3,7,11,13,p,q$ s.t.\ ${\ell} \nmid \nu_q(s - 1)$ (in particular: any prime $\ell \gg 0$), the absolute Galois group $G_k$ acts as a Zariski dense subgroup of $G_2(\bar{\mathbb{Q}}_{\ell})$ on the $\ell$-adic realization $H_\ell^\bullet(M_s)$ of the motive for motivated cycles $M_s$.
\end{thm}

\proof
By Cor.\ \ref{Cor: motivic interpretation}, we have to show that the image of the specialization $\rho_{\mathcal{V}_k(3)}^{(s)}$ at $s$ is a Zariski dense subgroup of $G_2(\bar{\mathbb{Q}}_{\ell})$.
The images of $\rho_{\mathcal{V}_k(3)}$ and $\rho_{\mathcal{V}_{\mathbb{Q}(\zeta_7)}^\varepsilon}$ lie inside $G_2(\bar{\mathbb{Q}}_{\ell})$ by Lem.\ \ref{lem: choice of the quadratic extension} and the definition of the character $\varepsilon$.
Since $G_2$ is a simple group, it suffices to show the analogue claim of the proposition for $\rho_{\mathcal{V}_{\mathbb{Q}(\zeta_7)}^\varepsilon}$. 
\\
By Lem.\ \ref{lem: irreducibility of the tate realization}, $G_{\mathbb{Q}(\zeta_7)}$ acts as an irreducible subgroup of $G_2(\bar{\mathbb{Q}}_{\ell})$ on $\bar{s}^*\mathcal{V}_{\mathbb{Q}(\zeta_7)}^{\varepsilon}$.
Consider the classification of the Zariski closed maximal subgroups $M$ of $G_2(\bar{\mathbb{Q}}_\ell)$ in \cite{Aschbacher} Cor.\ 12:
Say $M$ contains the image of $\rho_{\mathcal{V}_{\mathbb{Q}(\zeta_7)}^{\varepsilon}}^{(s)}$.
If $M$ is of type (1) to (3) in this classification, then it is reducible, which contradicts Lem.\ \ref{lem: irreducibility of the tate realization}.
Type (4) could not happen, since we work in characteristic $0$.
Thus, only type (5) remains, i.e., $M \cong {\rm PSL}_2(\bar{\mathbb{Q}}_\ell)$ and the action on $s^*\mathcal{V}_{\mathbb{Q}(\zeta_7)}^{\varepsilon}$ is isomorphic to the ${\rm PSL}_2(\bar{\mathbb{Q}}_\ell)$-action on the homogeneous polynomials of degree $6$ in $\bar{\mathbb{Q}}_\ell[X,Y]$ given as in (\ref{eq: PSL-action}).
We will show that this is not possible, too:
\\
The Jordan canonical form of every non trivial unipotent element in the image of this ${\rm PSL}_2(\bar{\mathbb{Q}}_\ell)$-representation consists of a single Jordan block of length $7$.
Arguing as in the proof of Lem.\ \ref{lem: specialization}, we see that the image of $\rho_{\mathcal{V}_{\mathbb{Q}(\zeta_7)}^{\varepsilon}}^{(s)}$ contains $g_1^{\nu_q(s-1)}$.
The Jordan canonical form of $g_1$ consists of three Jordan blocks, i.e., $g_1$ has a $3$-dimensional $1$-eigenspace.
It follows that $g_1^{\nu_q(s-1)}$ has at least a $3$-dimensional $1$-eigenspace, i.e., its Jordan canonical form consists of at least three blocks.
But $g_1^{\nu_q(s-1)}$ is non trivial unipotent, so its Jordan canonical from should consist only of a single block, a contradiction.
\\
Thus, no maximal Zariski closed subgroup of $G_2(\bar{\mathbb{Q}}_\ell)$ contains the image of $\rho_{\mathcal{V}_{\mathbb{Q}(\zeta_7)}^{\varepsilon}}^{(s)}$, i.e., the closure of this image is $G_2(\bar{\mathbb{Q}}_\ell)$ itself, which finishes the proof.
\Endproof

Suppose the Tate conjecture holds.
In particular, the standard conjectures on algebraic cycles hold and $\underline{\rm NM}(k)_{E}$ is Tannakian with fibre functor $H_{\sigma}^\bullet$.
Recall that the motivic Galois group $G_{M_s}^{\rm mot}$ is the Tannaka group of the $\otimes$-subcategory generated by $M_s$ (with resp.\ to the restriction of the Betti realization $H_\sigma^\bullet$).
It is a Zariski closed subgroup of ${\rm GL}(H_\sigma^\bullet(M_s))$ by \cite{DelMil82} Prop.\ 2.20 (b).
It is even reductive by Jannsen's Theorem (\cite{Jannsen92} Thm.\ 1) together with \cite{DelMil82} Prop.\ 2.23.
Thus, the base extension $G_{M_s}^{\rm mot} \otimes_E \bar{\mathbb{Q}}_\ell$ is exactly given as the fixgroup of all of its invariants in the induced representations $H_\ell^\bullet(M_s)^{\otimes n} \otimes H_\ell^\bullet(M_s)^{\vee,\otimes m}$ for $n,m \geq 0$
(cf.\ \cite{Andre04}).
By the Tate conjecture, these are exactly all the Tate-cycles in $H_\ell^\bullet(M_s)^{\otimes n} \otimes H_\ell^\bullet(M_s)^{\vee,\otimes m}$.

\begin{cor}\label{cor: numerical motives under tate conjecture}
Suppose the Tate conjecture holds.
Let $s \in \mathbb{Q} \setminus \{ 0,1 \}$ and $k/\mathbb{Q}(\zeta_7)$ be as in Thm.\ \ref{thm: tate realization}. 
Then $M_s$ is a well defined numerical motive over $k$ with motivic Galois group of type $G_2$.
\end{cor}

\proof
Choose a prime $\ell$ for $s$ as in Thm.\ \ref{thm: tate realization}.
It suffices to show that $G_{M_s}^{\rm mot}(\bar{\mathbb{Q}}_\ell)$ is isomorphic to $G_2(\bar{\mathbb{Q}}_{\ell})$.
The $G_k$-representation on the $\ell$-adic realization $H_{\ell}^\bullet(M_s) = s^*\mathcal{V}_k$ has image Zariski dense inside $G_2(\bar{\mathbb{Q}}_\ell)$.
Since this representation fixes all Tate-cycles on all $\otimes$-constructions of $H_{\ell}^\bullet(M_s)$ by definition, $G_{M_s}^{\rm mot}(\bar{\mathbb{Q}}_\ell)$ contains $G_2(\bar{\mathbb{Q}}_\ell)$ by the remarks preceding the corollary.
On the other hand, it fixes the underlying bilinear form $b \in S^2 H_{\ell}^\bullet(M_s)^\vee$ and the Dickson form $f \in \Lambda^3 H_{\ell}^\bullet(M_s)^\vee$, i.e.\ $b$ and $f$ are Tate cycles.
It follows that $G_{M_s}^{\rm mot}(\bar{\mathbb{Q}}_\ell)$ fixes these cycles, too, hence is contained inside $G_2(\bar{\mathbb{Q}}_{\ell})$ by \cite{Aschbacher} (2.11) and (3.4).
\Endproof

\begin{rem} 
In the last result, the underlying motive can be defined over $\mathbb{Q}(\zeta_7)$ using a suitable twist with an Artin motive 
of degree $\leq 2,$ given by  Lem.~\ref{lem: choice of the quadratic extension}.
\end{rem}

\subsection{An unconditional result for motives for absolute Hodge-cycles.}\label{sect: unconditional result}
In this final section, we want to sketch the proof for an unconditional analogue of Cor.\ \ref{cor: numerical motives under tate conjecture} for motives for absolute Hodge-cycles.
Now that we understand the Galois-action on the stalk $\bar{s}^*\mathcal{V}_k(3)$ (see Thm.\ \ref{thm: tate realization}), the proof itself goes almost literally as the arguments leading to \cite{Dettweiler1} Thm.\ 3.3.1, so we just deal with the differences.
The key difference is that we have to work with coefficients $E= \mathbb{Q}(\zeta_7)$ over the base field ${K}=\mathbb{Q}(\zeta_7)$.
Since $\mathbb{Q}$ is the base field in \cite{Dettweiler1}, ``$7^{\rm th}$-root'' of the determinant motive occurring in the proof of \cite{Dettweiler1} Thm.\ 3.3.1 could be trivialized by a Tate-twist using \cite{KisinWortmann} Thm.\ 3.1.
Unfortunately, the method of \cite{KisinWortmann} does not work over $K$ and coefficients $E$:
A priori it might happen that the ${K}\otimes E$-linear filtration of the de Rham realization of a rank one motive has filtration steps given by non trivial idempotents of $K\otimes E$.
So in our case, we have to show ``by hand'' that this does not happen for the determinant motive mentioned above.

\begin{thm}\label{thm: motives for absolute hodge cycles}
Let $s\in \mathbb{Q}\setminus \{ 0,1 \}$ be a rational number and suppose there are prime numbers $p,q \neq 2, 7$ s.t.\ $p \equiv 3 ~{\rm or}~ 5 ~{\rm mod}~ 7$, $\nu_p(s) < 0 < \nu_q(s - 1)$ and $7 \nmid \nu_p(s)$.
Let $k/\mathbb{Q}(\zeta_7)$ be the extension of degree $\leq 2$ given in Lem \ref{lem: choice of the quadratic extension}.
We identify $s$ with the induced $k$-point of $\mathbb{A}^1\setminus \{ 0,1 \}$.
\begin{enumerate}
 \item As a motive for motivated cycles over $k$ and up to a twist by a suitable rank one motive, $M_s$ has motivic Galois group of type $G_2$.
 \item As a motive for absolute Hodge-cycles over $k$, $M_s$ itself has motivic Galois group of type $G_2$.
\end{enumerate}
\end{thm}

\begin{rem}
The rank one motive in Thm.\ \ref{thm: motives for absolute hodge cycles} is explicitly given via the Tannaka formalism of the $\otimes$-category generated by $M_s$.
\end{rem}

\proof
For $M$ a motive for motivated resp.\ absolute Hodge-cycles we denote by $G_M^{\rm mot}$ resp.\ $G_M^{\rm AH}$ the (unconditional) motivic Galois group.
The cycle class of a motivated resp.\ absolute Hodge-cycle is by definition a Tate cycle.
Fix a $k$-point as in the assumption of the theorem and choose a prime $\ell$ as in Thm.\ \ref{thm: tate realization}.
As in the proof of Cor.\ \ref{cor: numerical motives under tate conjecture} it follows that both $G_{M_s}^{\rm mot}(\bar{\mathbb{Q}}_\ell)$ and $G_{M_s}^{\rm AH}(\bar{\mathbb{Q}}_\ell)$ contain $G_2(\bar{\mathbb{Q}}_\ell)$.
\\
Using \cite{Andre96} Thm.\ 5.2 literally as in the proof of \cite{Dettweiler1} Thm.\ 3.3.1, we get that the motive for motivated cycles $N_r$ over $\mathbb{C}$ has motivic Galois group
\begin{equation*}
 G_{N_r}^{\rm mot}(\bar{\mathbb{Q}}_\ell) \leq G_2(\bar{\mathbb{Q}}_\ell) \times {\rm scalar~matrices}
\end{equation*}
for any $r \in \mathbb{A}^1(\mathbb{C})\setminus \{0,1\}$ (we apply \cite{Andre96} Thm.\ 5.2 to the family $\{N_r\}_r$ seen as a family of rational motives with an $E$-module structure - this suffices for our needs since the $E$-linear motivic Galois group $G_{N_r}^{\rm mot}$ is contained in the $\mathbb{Q}$-linear one).
The same follows for the third Tate twist $M_r$.
Let $A_r$ be the motive dual to the motive corresponding to the resulting character $G_{N_r}^{\rm mot} \rightarrow \mathbb{G}_{\rm m}$ under the Tannaka formalism.
Note that $A_r^{\otimes 7}$ is by construction dual to the the determinant motive ${\rm det}M_r = \Lambda^7M_r$.
Further, $M_r \otimes A_r$ has motivic Galois group (with resp.\ to motivated cycles) contained inside $G_2$.
In particular, $b \in S^2 H_{\ell}^\bullet(M_r \otimes A_r)^\vee$ and the Dickson form $f\in\Lambda^3 H_{\ell}^\bullet(M_r \otimes A_r)^\vee$ are motivated cycles and hence absolute Hodge-cycles over $\mathbb{C}$, hence also over $\bar{k}$ (see the Scolie in Sect.\ 2.5.\ of \cite{Andre96} resp.\ \cite{Del82} Prop.\ 2.9).
By construction, the $\ell$-adic realization of $M_s$ is just $\bar{s}^*\mathcal{V}_k$.
It follows from Thm.\ \ref{thm: tate realization} that $b$ and $f$ are also Tate-cycles over $k$.
Thus, again by the Scolie in Sect.\ 2.5.\ of \cite{Andre96}, $b$ and $f$ are motivated cycles and hence absolute Hodge-cycles over $k$ and $G_{M_s \otimes A_s}^{\rm mot}(\bar{\mathbb{Q}}_\ell)$ resp.\ $G_{M_s \otimes A_s}^{\rm AH}(\bar{\mathbb{Q}}_\ell)$ is contained in $G_2(\bar{\mathbb{Q}}_\ell)$.
The above character is trivial on $G_2(\bar{\mathbb{Q}}_\ell)$ and $G_{M_s}^{\rm mot}(\bar{\mathbb{Q}}_\ell)$ contains $G_2(\bar{\mathbb{Q}}_\ell)$ by the above, so $G_{M_s \otimes A_s}^{\rm mot}(\bar{\mathbb{Q}}_\ell)$ is in fact $G_2(\bar{\mathbb{Q}}_\ell)$ and we are done with the first claim.
\\
For the second claim, it suffices to show that $A_s$ is trivial as a motive for absolute Hodge-cycles over $k$.
To do this, it suffices to see that the Galois action on the $\lambda$-adic realizations for rational primes $\lambda$ as well as the filtration on the de Rham realization is trivial.
\\
We start with the $\lambda$-adic realizations:
By construction these are of weight $0$ so it suffices to see that the characteristic polynomials of the Frobenius-actions are trivial.
Let $v$ be a place of $R_N$ not lying over $\lambda$.
The construction of $\mathcal{V}_{R_N}$ for various $\lambda$ in \cite{Dettweiler1} Thm.\ 2.4.1 is independent from $\lambda$, so the same is true for the characteristic polynomial of the Frobenius $F_v$ acting on $\bar{s}^*\mathcal{V}_k \cong H_\lambda^\bullet(M_s)$ and hence also for the Frobenius-action on the $\otimes$-construction $H_\lambda^\bullet(A_r)$.
Thus, it suffices to treat the case $\lambda = \ell$.
In this case, the Galois action on $H_\ell^\bullet(A_s)$ is given by the character $\varepsilon(-3)$, i.e., is trivial by our choice of $k$ according to Lem.\ \ref{lem: choice of the quadratic extension}. 

For the de Rham realization  we argue as follows:
If $e$ is a non trivial idempotent of $E\otimes K$ cutting out a non trivial filtration step of $H_{\rm dR}^\bullet(A_s)$, then $e=e^7$ cuts out a non trivial filtration step in the filtration of $H_{\rm dR}^\bullet({\rm det}M_s)^\vee = H_{\rm dR}^\bullet(A_s)^{\otimes 7}$.
Thus, we are done if we can show that $H_{\rm dR}^\bullet({\rm det}M_s)$ carries the trivial filtration.
\\
Everything is already defined over ${K}$, so we may work over the smaller field ${K} \leq k$.
By Rem.\ \ref{rem: motive defined over the rationals} we have
\begin{align*}
 F^qH_{\rm dR}^\bullet(M_s) &= {\rm Eig}(\omega_{\vert F^qH_{\rm dR}^\bullet(M_s^0) \otimes {K}} \otimes E, 1 \otimes \zeta_{14})
\\
 &\cong \bigoplus_{g \in {\rm Gal}(E/\mathbb{Q})} {\rm Eig}(\omega_{\vert F^qH_{\rm dR}^\bullet(M_s^0) \otimes {K}} , \zeta_{14}^g),
\end{align*}
where we treat $M_s^0$ as a rational motive over $\mathbb{Q}$.
It follows that
\begin{equation*}
 H_{\rm dR}^\bullet({\rm det}M_s) \cong \bigoplus_{g \in {\rm Gal}(E/\mathbb{Q})} \Lambda_{K}^7 {\rm Eig}(\omega_{\vert H_{\rm dR}^\bullet(M_s^0) \otimes {K}} , \zeta_{14}^g) 
\end{equation*}
holds as a filtered $K\otimes E$-module, where the filtration on the right hand side is the direct sum filtration.
Now the filtration on the direct summands is induced by the filtration on $H_{\rm dR}^\bullet(M_s^0)$, which is trivial since it corresponds to a rational Hodge structure  of weight $0$.
\Endproof

\bibliographystyle{alpha}

\begin{thebibliography}{DRK10}

\bibitem[And96]{Andre96}
Yves Andr\'{e}.
\newblock {Pour une th\'{e}orie inconditionelle des motifs}.
\newblock {\em Publ. Math. IHES}, 83:\ 5--49, 1996.

\bibitem[And04]{Andre04}
Yves Andr\'{e}.
\newblock {\em Une Introduction aux Motifs (Motifs Purs, Motifs Mixtes,
  P\'{e}riodes)}, volume~17 of {\em Panoramas et synth\`{e}ses}.
\newblock Soci\'{e}t\'{e} Math\'{e}matique de France, 2004.

\bibitem[Asc87]{Aschbacher}
Michael Aschbacher.
\newblock {Chevalley Groups of Type $G_2$ as the Group of a Trilinear Form}.
\newblock {\em Journal of Algebra}, 109:\ 193--259, 1987.

\bibitem[Bec91]{Beckmann}
Sybilla Beckmann.
\newblock {On extensions of number fields obtained by specializing branched
  coverings}.
\newblock {\em Journal f\"{u}r die reine und angewandte Mathematik}, 419:\
  27--53, 1991.


\bibitem[Car72]{Carter}
Roger~W. Carter.
\newblock {\em Simple Groups of Lie Type}.
\newblock Pure and Applied Mathematics. John Wiley \& Sons, London-New
  York-Sydney-Toronto, 1972.

\bibitem[Car85]{Carter2}
Roger~W. Carter.
\newblock {\em Finite Groups of Lie Type}.
\newblock Pure and Applied Mathematics; A Wiley-Interscience publication. John
  Wiley \& Sons, Chichester [u.a.], 1985.
  
\bibitem[CLS11]{CoxLittleSchenck}
D.~Cox, J.~Little and H.~Schenck.
\newblock {\em Toric Varieties}, volume~124 of {\em Graduate Studies in Mathematics}.
\newblock American Mathematical Society, 2011.

\bibitem[CR74]{changree}
Bomshik Chang and Rimhak Ree.
\newblock {The Characters of $G_2(q)$}.
\newblock {\em Symposia Mathematica}, XIII:\ 395--413, 1974.

\bibitem[Del80]{WeilII}
Pierre Deligne.
\newblock {La conjecture de Weil. II.}
\newblock {\em Publ. Math. IHES}, 52:\ 137--252, 1980.

\bibitem[Del82]{Del82}
Pierre Deligne (Notes by James S.~Milne).
\newblock {Hodge cycles on abelian varieties}.
\newblock In {\em Hodge cycles, motives, and Shimura varieties}, volume
  900 of {\em Lecture Notes in Mathematics}, pages 9--100. Springer-Verlag,
  Berlin-Heidelberg-New York, 1982.

\bibitem[DM82]{DelMil82}
Pierre Deligne and James S.~Milne.
\newblock {Tannakian Categories}.
\newblock In {\em Hodge cycles, motives, and Shimura varieties}, volume
  900 of {\em Lecture Notes in Mathematics}, pages 101--228. Springer-Verlag,
  Berlin-Heidelberg-New York, 1982.

\bibitem[DR00]{Dettweiler3}
Michael Dettweiler and Stefan Reiter.
\newblock {An Algorithm of Katz and its Application to the Inverse Galois
  Problem}.
\newblock {\em Journal of Symbolic Computation}, 30:\ 761--798, 2000.

\bibitem[DRK10]{Dettweiler1}
Michael Dettweiler, Stefan Reiter, and Nicholas~M.\ Katz.
\newblock {Rigid local systems and motives of type $G_2$. With an appendix by
  Michael Dettweiler and Nicholas M.\ Katz}.
\newblock {\em Compositio Mathematica}, 146(4):929--963, 2010.

\bibitem[ea85]{Atlas}
J.H.~Conway et~al.
\newblock {\em Atlas of finite groups}.
\newblock Clarendon Press, Oxford, 1985.

\bibitem[ESV92]{EsnaultSchechtmanViehweg}
H.~Esnault, V.~Schechtman and E.~Viehweg.
\newblock {Cohomology of local systems on the complement of hyperplanes}.
\newblock {\em Inventiones Mathematicae}, 109(1):\ 557--5, 1992.

\bibitem[Jan92]{Jannsen92}
Uwe Jannsen.
\newblock {Motives, numerical equivalence, and semi-simplicity}.
\newblock {\em Inventiones Mathematicae}, 107(1):\ 447--452, 1992.

\bibitem[Kat96]{Katz}
Nicholas~M. Katz.
\newblock {\em Rigid Local Systems}, volume 139 of {\em Annals of Mathematics
  Studies}.
\newblock Princeton University Press, Princeton, New Jersey, 1996.

\bibitem[KW03]{KisinWortmann}
M.~Kisin and S.~Wortmann.
\newblock {A note on Artin motives}.
\newblock {\em Mathematical Research Letters}, 10(3):\ 375--389, 2003.

\bibitem[Kle88]{Kleidman}
Peter~B. Kleidman.
\newblock {The Maximal Subgroups of the Chevalley Groups $G_2(q)$ with $q$ Odd,
  the Ree Groups $^2 G_2(q)$, and Their Automorphism Groups}.
\newblock {\em Journal of Algebra}, 117:\ 30--71, 1988.

\bibitem[Mal03]{Malle}
Gunter Malle.
\newblock {Explicit realization of the Dickson groups $G_2(q)$ as Galois
  groups}.
\newblock {\em Pacific Journal of Mathematics}, 212(1):\ 157--167, 2003.

\bibitem[MM99]{Matzat}
G.~Malle and B.H. Matzat.
\newblock {\em Inverse Galois Theory}.
\newblock Springer Monographs in Mathematics. Springer-Verlag,
  Berlin-Heidelberg-New York, 1999.


\bibitem[Ser94]{Serre}
J.-P. Serre.
\newblock {Propri\'{e}t\'{e}s conjecturales des groupes de Galois motiviques et
  des repr\'{e}sentations $\ell$-adiques}.
\newblock {\em Proceedings of Symposia in Pure Mathematics}, 55:\ 377--400,
  1994.
  
\bibitem[SGA1]{SGAI}
Alexander~Grothendieck et~al.
\newblock {\em S\'{e}minaire de g\'{e}om\'{e}trie alg\'{e}brique du Bois-Marie
  1960-61, Rev\^{e}tements Etales et Group Fondamental (SGA 1)}, volume 224 of
  {\em Lecture Notes in Mathematics}.
\newblock Springer-Verlag, Berlin-Heidelberg-New York, 1971.


\bibitem[Yun14]{Yun14}
Zhiwei Yun.
\newblock {Motives with exceptional Galois groups and the inverse Galois problem}.
\newblock {\em Inventiones Math.}, 196(2):\ 267-337, 2014.

\end{thebibliography}

\end{document}